\documentclass[reqno,10pt]{amsart}
\usepackage{amsfonts}
\usepackage{graphicx}
\usepackage{amsfonts,amsmath, amssymb}
\usepackage{hyperref}
\usepackage{color}
\usepackage{cite}
\usepackage{bm}
\usepackage{enumerate}

\def\k{\kappa}

\def\a{\alpha}
\def\g{\gamma}

\def\s{\sigma}

\def\f{\frac}

\numberwithin{equation}{section}

\theoremstyle{plain}
\newtheorem{theorem}{Theorem}[section]

\newtheorem{lemma}[theorem]{Lemma}
\newtheorem{proposition}[theorem]{Proposition}

\theoremstyle{remark}
\newtheorem{remark}[theorem]{Remark}

\def\k{\kappa}

\def\a{\alpha}
\def\g{\gamma}

\def\s{\sigma}

\def\f{\frac}

\renewcommand{\P}{\mathbb{P}}

\renewcommand{\S}{\mathbb{S}}

\newcommand{\Fi}{\mathbf{1}}

\newcommand{\CF}{\mathcal{F}}

\newcommand{\CL}{\mathcal{L}}

\newcommand{\CM}{\mathcal{M}}

\newcommand{\CQ}{\mathcal{Q}}
\newcommand{\CV}{\mathcal{V}}

\newcommand{\MRB}{\mathfrak{B}}

\newcommand{\ga}{\gamma}

\newcommand{\pa}{\partial}
\newcommand{\ka}{\kappa}

\newcommand{\vep}{\varepsilon}

\newcommand{\dd}{{\rm d}}
\newcommand{\eqdef}{\overset{\mbox{\tiny{def}}}{=}}

\begin{document}
\title[Boltzmann equation in a channel]{ Heat transfer problem  for the Boltzmann equation in a channel with diffusive boundary condition}

\author[R.J. Duan]{Renjun Duan}
\address[R.J. Duan]{Department of Mathematics, The Chinese University of Hong Kong, Shatin, Hong Kong, P. R. China}
\email{rjduan@math.cuhk.edu.hk}

\author[S.Q. Liu]{Shuangqian Liu}
\address[S.Q. Liu]{School of Mathematics and Statistics and Hubei Key Laboratory of Mathematical Sciences, Central China Normal University, Wuhan 430079, P. R. China}
\email{tsqliu@jnu.edu.cn}
\author[T. Yang]{Tong Yang}
\address[T. Yang]{Department of Mathematics, City University of Hong Kong, Hong Kong, P. R. China}
\email{matyang@cityu.edu.hk}

\author[Z. Zhang]{Zhu Zhang}
\address[Z. Zhang]{Department of Mathematics, City University of Hong Kong, Hong Kong, P. R. China}
\email{zhuzhangpde@gmail.com}
\maketitle

\begin{center}{\large \it
		Dedicated to the memory of Professor Chaohao Gu}
\end{center}
\begin{abstract}
In this paper, we study the 1D steady Boltzmann flow in a channel. The walls of the channel are assumed to have  vanishing velocity and  given temperatures $\theta_0$ and $\theta_1$. This problem was studied by Esposito et al \cite{ELM,ELM2} where they showed that the solution tends to a local Maxwellian with parameters satisfying the compressible Navier-Stokes equation  with no-slip boundary condition. However, a lot of numerical experiments reveal that the fluid layer does not entirely stick to the boundary. In the regime where the Knudsen number is reasonably small, the slip phenomenon is significant near the boundary. Thus, we revisit this problem by taking into account the slip boundary conditions. Following the lines of \cite{Coron}, we will first give a formal asymptotic analysis to see that the flow governed by the Boltzmann equation  is accurately approximated by a superposition of a steady CNS equation with a temperature jump condition and two Knudsen layers located at end points. Then we will establish a uniform $L^\infty$ estimate on the remainder and derive the slip boundary condition for compressible Navier-Stokes equations rigorously.
\end{abstract}

\section{Introduction}
\subsection{Problem settings}
In this paper, we study the steady flow of a rarefied gas in a channel which is bounded by two thermal walls located at $x=0$ and $x=1.$ The walls are assumed to have a vanishing  velocity and  given temperatures $\theta_{0}$ and $\theta_1$ ($\theta_0\neq \theta_1$) respectively. In the kinetic setting, the distribution function satisfies the following 1D rescaled steady Boltzmann equation:
\begin{align}\label{C.1.0.1}
v_1\pa_xF^\vep=\vep^{-1}Q(F^\vep,F^\vep),\quad x\in (0,1),~ v=(v_1,v_2,v_3)\in \mathbb{R}^3.
\end{align}
The parameter $\vep>0$ is the Knudsen number which is proportional to the mean free path and is assumed to be small.
The Boltzmann collision term on the right-hand side of \eqref{C.1.0.1} takes the non-symmetric bilinear form of
\begin{equation}
Q(F_1,F_2)=\int_{\mathbb{R}^3}\int_{\mathbb{S}^2} B(|v-u|,\omega)[F_1(u')F_2(v')-F_1(u)F_2(v)]\,\dd\omega\dd u,\nonumber
\end{equation}
where the velocity pair $(v',u')$ is defined by the velocity pair $(v,u)$ as well as the parameter $\omega\in \S^2$ in  the relation
\begin{equation*}
v'=v-[(v-u)\cdot\omega]\omega,\quad u'=u+[(v-u)\cdot\omega]\omega,
\end{equation*}
according to conservation laws of momentum and energy
\begin{equation*}
v'+u'=v+u,\quad |v'|^2+|u'|^2=|v|^2+|u|^2
\end{equation*}
for elastic collision. For simplicity, we consider the hard sphere model when the collision kernel $B(|v-u|,\omega)=|v-u||\cos\psi|$ with  relative velocity $v-u$ and interaction angle defined by
$\cos \psi=\omega \cdot (v-u)/|v-u|$. Denote the outward normal vectors at each boundary point $x\in \{0,1\}$ by
$$n(x)=\left\{\begin{aligned}
(-1,0,0),\quad &x=0,\\
(1,0,0),\quad &x=1.
\end{aligned}\right.
$$
The phase boundary $\{0,1\}\times \mathbb{R}^3$ can be decomposed into three parts:
$$\{0,1\}\times \mathbb{R}^3=\g_+\cup \g_0\cup\g_-,$$ where
\begin{align}
\g_{\pm}&=\left\{(x,v)\mid n(x)\cdot v\gtrless 0\right\}\nonumber\\
&=\left(\{0\}\times \{v_1\lessgtr0\}\right)\cup\left(\{1\}\times \{v_1\gtrless0\}\right),
\label{C.a.10}\end{align} and $$\g_0=\{0,1\}\times \{v_1=0\}.$$
The distribution of gas particles at thermal walls satisfies the diffusive reflection boundary condition, which is given by
\begin{equation}
\label{C.1.0.2}
F^\vep(x,v)|\g_-=\mu_{\theta_w}\int_{\{n(x)\cdot v>0\}}F^\vep(x,v)\{n(x)\cdot v\}\,\dd v,\quad x=0,1,
\end{equation}
where the wall Maxwellian $\mu_{\theta_w}$ is 
\begin{align}\nonumber
\mu_{\theta_{w}}(v)=\f{1}{2\pi\theta^2_w(x)}e^{-\f{|v|^2}{2\theta_w(x)}},~\theta_{w}(0)=\theta_0,~ \theta_w(1)=\theta_1.
\end{align}
Without loss of generality, we assume that $\theta_0=1$ and $\theta_1>\theta_0.$  Moreover, the total mass of the solution is  equal to $1$ throughout the paper, i.e.
\begin{align}
	\label{C.2.0.2}
	\int_0^1\int_{\mathbb{R}^3}F^\vep(x,v)\,\dd v\dd x=1.
\end{align}
\subsection{CNS approximation}
We are interested in the behavior of solution $F^\vep$ in the limit $\vep \rightarrow 0^+$ that  is the hydrodynamic limit of the Boltzmann equation. In the absence of physical boundaries or shocks, it is well-known that the distribution function converges to a local Maxwellian with parameters satisfying the compressible Euler system; {\it cf.} \cite{Gra2}. The Chapman-Enskog expansion yields the compressible Navier-Stokes system (CNS) as the first order correction. In this subsection, we give a formal derivation of CNS approximation in the setting of this paper. Before this, we  define some function spaces which will be used later. Given a local Maxwellian $$\CM=\CM_{[\rho,u,\theta]}(v)=\f{\rho}{(2\pi\theta)^{3/2}}e^{-\f{|v-u|^2}{2\theta}}$$ with density $\rho(x)>0$, velocity $ u(x)$ and temperature $\theta(x)>0$, we define a function space $L^2_{\CM}$, equipped with the following inner product in $v$:
$$\langle f,g\rangle_{\CM}=\int_{\mathbb{R}^3}\f{fg}{\CM}\,\dd v.
$$  The linearized collision operator around $\CM$ is given by
\begin{align}
\CL_\CM(\cdot)=-\left[Q(\CM,\cdot)+Q(\cdot,\CM)\right].\nonumber
\end{align} The following properties of $\CL_{\CM}$ are well-known (cf. \cite{Gra}). That is,  $\CL_{\CM}$ has a null space $\text{Ker}\CL_{\CM}$ in $L^2_{\CM}$ that  is spanned by the following functions
$$\chi_1=\CM,~\chi_i= \frac{v_i}{\sqrt{\theta}}\CM,~i=1,2,3,~\chi_4=\frac{|v|^2-3\theta}{\sqrt{6}\theta}\CM.
$$
We also define the macroscopic projection operator $\P_{\CM}$ as the projection onto $\text{Ker}\CL_{\CM}$. If  $f$ is orthogonal to null space of $\CL_\CM$, the following coercive inequality holds:
\begin{align}\label{C.a.12}
\langle \CL_\CM f, f\rangle_{L^2_\CM}\geq c_0|\nu^{\f12}f|_{L^2_\CM},
\end{align}
for some positive constants $c_0$, where $\nu(v)=1+|v|$.

Without boundary, by formally passing limit $\vep\rightarrow 0$ in \eqref{C.1.0.1} we see that the leading order gives the local Maxwellian $\CM$ which satisfies $Q(\CM,\CM)=0$. Hence,
we set the following expansion \begin{align}\label{C.1.1.1}
F^\vep= \CM+\vep G+\vep^2 F_2+O(\vep^3),
\end{align}
where $G$ and $F_2$ are some correctors to be constructed later.
Inserting the expansion \eqref{C.1.1.1} into \eqref{C.1.0.1} yields the error
$$\begin{aligned}E=v_1\pa_x\CM+\CL_MG+\vep\left[v_1\pa_xG+\CL_{\CM}F_2-Q(G,G)\right]+O(\vep^2).
\end{aligned}$$
Then we can eliminate the terms of $O(1)$ and $O(\vep)$ order by choosing
\begin{align}
G&=-\CL_{\CM}^{-1}\left[(I-\P_{\CM})v_1\pa_x\CM\right],\label{G}\\
F_2&=-\CL_{\CM}^{-1}\left[(I-\P_{\CM})v_1\pa_xG-Q(G,G) \right]+\tilde{F}_2,\label{F2}
\end{align}
where $\tilde{F}_2\in \text{Ker}\CL_\CM$ will be determined later. Here, 
$\CM$ satisfies
\begin{align}\label{C.1.1.2}
\P_{\CM}v_{1}\pa_x\CM+\vep \P_{\CM}v_1\pa_xG=0.
\end{align}
It is straightforward to check that \eqref{C.1.1.2} is equivalent to the following 1-D steady compressible Navier-Stokes system for $x\in (0,1)$:
\begin{equation}\label{C.1.1.3}
\left\{
\begin{aligned}
&\pa_{x}(\rho u_1)=0,\\
&\pa_{x}(\rho u_1^2+\rho\theta)=\vep\pa_x(\f{4}{3}\iota(\theta)\pa_{x}u_1),\\
&\pa_{x}(\rho u_1u_j)=\vep\pa_x(\iota(\theta)\pa_{x}u_j),\quad j=2,3,\\
&\pa_{x}\left[\rho u_1\left(\f{5\theta}{2}+\f{|u|^2}{2}\right)\right]=\vep\pa_x\left(\f{4}{3}\iota(\theta)u_1\pa_xu_1+\sum_{i=2,3}\iota(\theta)u_i\pa_xu_i+\ka(\theta)\pa_x\theta\right).
\end{aligned}\right.
\end{equation}
Here, the viscosity $\iota(\theta)$ and heat conductivity $\ka(\theta)$ are respectively given by
\begin{equation}
\left\{\begin{aligned}
&\iota(\theta)=\theta\int_{\mathbb{R}^3}A_j(\xi)\CL_{\CM_{[1,u,\theta]}}^{-1}A_j(\xi)\CM_{[1,u,\theta]}\dd v, \quad i=2 \text{ or }3,\\
&\ka(\theta)=\theta\int_{\mathbb{R}^3}B(\xi)\CL_{\CM_{[1,u,\theta]}}^{-1} B(\xi)\CM_{[1,u,\theta]}\dd v,~~\xi=\frac{v-u}{\sqrt{\theta}}, \nonumber
\end{aligned}\right.
\end{equation}
where $$A_j(\xi)=\xi_j\xi_1-\f{|\xi|^2}{3}\delta_{1j},\quad B(\xi)=\f{|\xi|^2-5}{2}\xi_1$$
are the Burnett functions. It is well-known (cf. \cite{BGL1}) that these functions satisfy the following relations
$$\begin{aligned}&-\CL_{\CM_{[1,u,\theta]}}^{-1}[A_{i}(\xi)\CM_{[1,u,\theta]}]=\alpha(|\xi|,\theta)A_i(\xi)\CM_{[1,u,\theta]},\\
&-\CL_{\CM_{[1,u,\theta]}}^{-1}[B(\xi)\CM_{[1,u,\theta]}]=\beta(|\xi|,\theta)B(\xi)\CM_{[1,u,\theta]},\end{aligned}
$$
with two smooth scalar functions $\alpha(|\xi|,\theta)$ and $\beta(|\xi|,\theta)$. After some
 scaling, one has $$\begin{aligned}&\alpha(\xi,\theta)A_i(\xi)\CM_{[1,u,\theta]}=\f{1}{\theta^{2}}\alpha(\xi,1)A_i(\xi)\mu,\\ &\beta(\xi,\theta)B(\xi)\CM_{[1,u,\theta]}=\f{1}{\theta^{2}}\beta(\xi,1)B(\xi)\mu,
\end{aligned}
$$
where
 $$\mu(\xi)=\f{1}{(2\pi)^{\f32}}e^{-|\xi|^2/2}$$ is the normalized global Maxwellian. Here,  the scalar functions
$\alpha(\xi,1), \beta(\xi,1)$ are determined by
$$\alpha(\xi,1)A_i(\xi)\mu=\CL_{\mu}^{-1} A_i(\xi)\mu,\quad \beta(\xi,1)B(\xi)\mu=\CL_{\mu}^{-1}B(\xi)\mu,
$$
and $\CL_\mu=-[Q(\mu,\cdot)+Q(\cdot,\mu)]$ is the linearized collision operator around $\mu$. Then viscosity and heat conductivity coefficients are given by
\begin{equation}
\left\{\begin{aligned}
&\iota(\theta)=\left(\int_{\mathbb{R}^3}|A_j(\xi)|^2\alpha(\xi,1)\mu(\xi)\dd \xi\right)\theta^{\f12}=\bar{\iota}\theta^{\f12}, \quad j=2 \text{ or }3,\\
&\ka(\theta)=\left(\int_{\mathbb{R}^3}|B(\xi)|^2\beta(|\xi|,1)\mu(\xi)\dd \xi\right)\theta^{\f12}=\bar{\ka}\theta^{\f12}. \nonumber
\end{aligned}\right.
\end{equation}
Similarly, $G$ can be rewritten as
\begin{align}\label{C.1.1.4}
G&=-\CL_{\CM_{[1,u,\theta]}}^{-1}\left[\sum_{i=1,2,3}A_{i}(\xi)\pa_xu_i+\f{\pa_x\theta}{\sqrt{\theta}}B(\xi)\right]\CM_{[1,u,\theta]}\nonumber\\
&=-\left[\alpha(|\xi|,\theta)\left(\sum_{i=1,2,3}A_i(\xi)\pa_xu_i\right)+\beta(|\xi|,\theta)B(\xi)\f{\pa_x\theta}{\sqrt{\theta}}\right]\CM_{[1,u,\theta]}\nonumber\\
&=-\f{1}{\theta^2}\left[\alpha(\xi,1)\left(\sum_{i=1,2,3}A_i(\xi)\pa_{x}{u}_i\right)+\beta(\xi,1)B(\xi)\f{\pa_x{\theta}}{\sqrt{\theta}}\right]\mu(\xi).
\end{align}

\subsection{Slip boundary condition}\label{sp-bcd}
In order to solve compressible Navier-Stokes system \eqref{C.1.1.3} when $x\in (0,1)$, suitable boundary conditions are needed.  If we consider the no-slip boundary condition $$u(0)=u(1)=0, \quad \theta(0)=\theta_0,\quad \theta(1)=\theta_1,$$ the approximation \eqref{C.1.1.1} matches the boundary conditions \eqref{C.1.0.2} up to $O(1)$. However, since $G$ contains non-Maxwellian terms, the Chapman-Enskog  approximation $\CM+\vep G$ in general does not match the boundary condition \eqref{C.1.0.2} up to $O(\vep)$, except for the case when
$$\pa_{x}u(0)=\pa_xu(1)=\pa_x\theta(0)=\pa_{x}\theta(1)=0.
$$ 
However, then \eqref{C.1.1.3} is overdetermined. To obtain a more accurate approximation, Coron \cite{Coron} {\it formally} derived the  slip boundary conditions for compressible Navier-Stokes equations, which are essentially   a consequence of the analysis of the Knudsen layer. In what follows, we elaborate the derivation only in one dimensional case. We refer to \cite{ABHKM,Sone2,Sone,SBGS} for
the  physical investigations in general cases.

As in \cite{Coron}, since Chapman-Enskog expansion is not valid near the boundary,
we introduce Knudsen layers $\mathfrak{B}_0$ and $\mathfrak{B}_1$ around boundary points $x=0$ and $x=1$ respectively. The construction of Knudsen layers relies on the solutions to the following Milne problem:
\begin{equation}\label{C.a.2}
\left\{\begin{aligned}
&v_1\pa_y\mathcal{F}+\CL_{\mu}\mathcal{F}=0, ~y>0,~v\in \mathbb{R}^3,\\
&\mathcal{F}(0,v)|_{v_1>0}=\mathcal{G},\\
&\int_{\mathbb{R}^3}v_1\mathcal{F}\dd v=0,\\
&\lim_{y\rightarrow\infty}\mathcal{F}(y)=\mathcal{F}_{\infty} \text{ exists and belongs to }  \text{Ker}\CL_{\mu},
\end{aligned}
\right.
\end{equation}
where $\mathcal{G}$ is a given incoming  distribution function. The well-posedness of \eqref{C.a.2} has been shown in \cite{BCN}, and is summarized in Lemma \ref{lmC.a.1} for later use.

Now we construct the Knudsen layer $\mathfrak{B}_0$ and $\mathfrak{B}_1$ at the boundary points $x=0$ and $x=1$ respectively. Let $p\in \{0,1\}$ be a boundary point.
We set the boundary conditions of Navier-Stokes system as
\begin{align}
\label{C.1.2.1-1}
u(p)=\vep\tilde{u}(p),~~ \theta(p)-\theta_p=\vep \tilde{\theta}(p),
\end{align}
where $\tilde{u}(p)$ and $\tilde{\theta}(p)$ are corrections to be determined later.
Then we expand the boundary values of $\CM(x,v)|_{x=p}$ and $G(x,v)|_{x=p}$ as:
\begin{align}\label{C.1.2.2}
\CM(p)=&\sqrt{\f{\theta_p}{2\pi}}\rho(p)\mu_{\theta_p}+\vep\left\{\f{\tilde{u}(p)\cdot v}{\theta_p}+\f{|v|^2-3\theta_p}{2\theta^2_p}\tilde{\theta}(p)\right\}\rho(p)\sqrt{\f{\theta_p}{2\pi}}\mu_{\theta_p}+\vep^2\CM_{R}(p)\nonumber\\
:=&\sqrt{\f{\theta_p}{2\pi}}\rho(p)\mu_{\theta_p}+\vep\CM_1(p)+\vep^2\CM_R(p),~~p=0,1,
\end{align}
and
\begin{align}
G(p)&=-\f{\mu_{\theta_p}}{\sqrt{2\pi}}\bigg\{\alpha\big(\f{|v|}{\sqrt{\theta_p}},1\big) A\big(\f{v}{\sqrt{\theta_p}}\big)\pa_x{u}(p)+\beta\big(\f{|v|}{\sqrt{\theta_p}},1\big) B\big(\f{v}{\sqrt{\theta_p}}\big)\f{\pa_x{\theta}(p)}{\sqrt{\theta_p}}\bigg\}+\vep G_R(p)\nonumber\\
:&
= G_0(p)+\vep G_R(p),~~p=0,1.\label{C.1.2.2-1}
\end{align}
Now we consider $\mathfrak{B}_0$ first.
To compensate $G_0$ at $x=0$, we take $\mathcal{F}_{\alpha,i},i=1,2,3$ and  $\mathcal{F}_{\beta}$ as  solutions to the Milne problem \eqref{C.a.2} with incoming distribution functions
 $\mathcal{G}_{\a,i}=\alpha(|v|,1)A_i(v)\mu, i=1,2,3$ and $\mathcal{G}_\beta=\beta(|v|,1)B(v)\mu$.
By $\eqref{C.a.2}_4$, there exist positive constants $c_{\a,1}, c_{\a,2}, \underline{c}, c_{\beta,1}$ and $c_{\beta,2}$, such that
\begin{equation}
\begin{aligned}
\mathcal{F}_{\alpha,1,\infty}&=\lim_{y\rightarrow +\infty}\mathcal{F}_{\alpha,1}(y)=c_{\alpha,1}\mu+c_{\alpha,2}\frac{|v|^2-3}{2}\mu,\\
\mathcal{F}_{\alpha,j,\infty}&=\lim_{y\rightarrow +\infty}\mathcal{F}_{\alpha,j}(y)=\underline{c}v_j\mu,~~j=2,3,\\
\mathcal{F}_{\beta,\infty}&=\lim_{y\rightarrow +\infty}\mathcal{F}_{\beta}(y)=c_{\beta,1}\mu+c_{\beta,2}\frac{|v|^2-3}{2}\mu.\nonumber
\end{aligned}
\end{equation}
Then we define the Knudsen layer $\mathfrak{B}_0$ as
\begin{align}
\mathfrak{B}_0=&\sum_{i=1,2,3}\left[\mathcal{F}_{\alpha,i}\left(\frac{\rho(0)x}{\vep},\frac{v}{\sqrt{\theta_0}}\right)-\mathcal{F}_{\alpha,i,\infty}\left(\frac{v}{\sqrt{\theta_0}}\right)\right]\frac{\pa_xu_i(0)}{\theta^2_0}\nonumber\\
&+\left[\mathcal{F}_{\beta}\left(\frac{\rho(0)x}{\vep},\frac{v}{\sqrt{\theta_0}}\right)-\mathcal{F}_{\beta,\infty}\left(\frac{v}{\sqrt{\theta_0}}\right)\right]\frac{\pa_x\theta(0)}{\theta^{5/2}_0}.\label{C.a.8}
\end{align}
By a straightforward computation, one has
\begin{align}\nonumber
	\mathfrak{B}_0|_{x=0}=-G_0(0)-\Psi(0),
\end{align}
where
	\begin{align}
	\label{C.1.2.5}\Psi(0)=&\f{c_{\alpha,1}\sqrt{\theta_0}\pa_{x}u_1(0)+c_{\beta,1}\pa_{x}\theta(0)}{\sqrt{2\pi\theta_0}}\mu_{\theta_0}+\sum_{i=2,3}\f{\underline{c}\pa_{x}{u}_i(0)v_i}{\sqrt{2\pi\theta_0}}\mu_{\theta_0}\nonumber\\
	&+\f{c_{\alpha,2}\sqrt{\theta_0}\pa_{x}{u_1(0)}+c_{\beta,2}\pa_{x}\theta(0)}{\sqrt{2\pi \theta_0}}\f{|v|^2-3\theta_0}{2\theta^2_0}\mu_{\theta_0}\in \text{Ker}L_{\mu_{\theta_0}}.
\end{align}

Note that the Maxwellian part in \eqref{C.1.2.5} already satisfies the boundary condition \eqref{C.1.0.2} at $x=0$. Then we use $\CM_1(0)$  given in \eqref{C.1.2.2} to eliminate non-Maxwellian terms in $\Psi(0)$. That is,  set
\begin{align}
\tilde{u}_1(0)=0,~~ \tilde{u}_i(0)=\f{\underline{c}\pa_xu_i(0)}{\rho(0)},i=2,3,~~\tilde{\theta}(0)=\f{c_{\alpha,2}\sqrt{\theta_0}\pa_{x}{u_1(0)}}{\rho(0)}+\f{c_{\beta,2}\pa_{x}\theta(0)}{\rho(0)}.\nonumber
\end{align}
In view of \eqref{C.1.2.1-1}, it requires $[u,\theta]$ to satisfy the following slip boundary condition:
\begin{align}\label{C.1.2.8}
&u_{1}(0)=0,\quad {\rho(0)}{u}_i(0)=\vep\underline{c}\pa_xu_i(0),i=2,3,\nonumber\\
& \rho(0)[\theta(0)-\theta_0]={\vep c_{\alpha,2}\sqrt{\theta_0}\pa_{x}{u_1(0)}}+{\vep c_{\beta,2}\pa_{x}\theta(0)}
\end{align}
at $x=0.$
Similarly, at $x=1$ we can construct Knudsen layer $\mathfrak{B}_1$ at $x=1$ as follows: 
\begin{align}
\mathfrak{B}_1=&\left[\mathcal{F}_{\alpha,1}\left(\frac{\rho(1)(1-x)}{\vep},\frac{\mathfrak{R}v}{\sqrt{\theta_1}}\right)-\mathcal{F}_{\alpha,1,\infty}\left(\frac{\mathfrak{R}v}{\sqrt{\theta_1}}\right)\right]\frac{\pa_xu(1)}{\theta^2_1}\nonumber\\
&+\sum_{j=2,3}\left[-\mathcal{F}_{\alpha,j}\left(\frac{\rho(1)(1-x)}{\vep},\frac{\mathfrak{R}v}{\sqrt{\theta_1}}\right)+\mathcal{F}_{\alpha,j,\infty}\left(\frac{\mathfrak{R}v}{\sqrt{\theta_1}}\right)\right]\frac{\pa_xu(1)}{\theta^2_1}\nonumber\\
&+\left[-\mathcal{F}_{\beta}\left(\frac{\rho(1)(1-x)}{\vep},\frac{\mathfrak{R}v}{\sqrt{\theta_1}}\right)+\mathcal{F}_{\beta,\infty}\left(\frac{\mathfrak{R}v}{\sqrt{\theta_1}}\right)\right]\frac{\pa_x\theta(1)}{\theta^{5/2}_1},\label{C.a.7}
\end{align}
where $\mathfrak{R}v=\mathfrak{R}(v_1,v_2,v_3)=(-v_1,v_2,v_3).$ Then we have
\begin{align}
\mathfrak{B}_1|_{x=1}=-G_0(1)-\Psi(1),\nonumber
\end{align}
where
\begin{align}\label{C.1.2.5-2}
\Psi(1)=&\f{-c_{\alpha,1}\sqrt{\theta_1}\pa_{x}u_1(1)+c_{\beta,1}\pa_{x}\theta(1)}{\sqrt{2\pi\theta_1}}\mu_{\theta_1}-\sum_{i=2,3}\f{\underline{c}\pa_{x}{u}_i(1)v_i}{\sqrt{2\pi\theta_1}}\mu_{\theta_1}\nonumber\\
&+\f{-c_{\alpha,2}\sqrt{\theta_1}\pa_{x}{u_1(1)}+c_{\beta,2}\pa_{x}\theta(1)}{\sqrt{2\pi \theta_1}}\f{|v|^2-3\theta_1}{2\theta^2_1}\mu_{\theta_1}\in \text{Ker}L_{\mu_{\theta_1}}.
\end{align}
Then we set the following boundary condition of $[u,\theta]$ at $x=1$:
\begin{align}\label{C.1.2.9}
&u_{1}(1)=0,\quad {\rho(1)}u_i(1)=-\vep\underline{c}\pa_xu_i(1),i=2,3,\nonumber\\ &\rho(1)[\theta_1-\theta(1)]=-{\vep c_{\alpha,2}\sqrt{\theta_1}\pa_{x}{u_1(1)}}+{\vep c_{\beta,2}\pa_{x}\theta(1)}.
\end{align}
It is straightforward to check that $\CM+\vep G+\vep \mathfrak{B}_1$ satisfies the boundary condition \eqref{C.1.0.2} at $x=1$, up to the order $\vep$.

\subsection{Main result}
The paper aims to justify rigorously the slip boundary conditions presented in the previous section. For this, we start with the following expansion
\begin{align}\label{C.1.3.1}
F^\vep=\CM+\vep G+\vep\mathfrak{B}_0+\vep\mathfrak{B}_1+\vep^2F_2+\vep^{1+\alpha}F_R,
\end{align}
where $\alpha>0$ is a positive constant. Here we elaborate the approximate solutions appearing in the expansion: The leading order term $\CM=\CM_{[\rho,u,\theta]}$ is a local Maxwellian where $[\rho,u,\theta]$ satisfies the steady compressible Navier-Stokes equations with {\it slip} boundary conditions \eqref{C.1.2.8} and \eqref{C.1.2.9}. It will be constructed in Sec. \ref{sec3.1}.  The function $G$ is a corrector at order $\vep$ which is defined in  \eqref{C.1.1.4} and it satisfies \eqref{G}.  $\mathfrak{B}_0$ and $\mathfrak{B}_1$ are Knudsen layers which are defined in  \eqref{C.a.8} and \eqref{C.a.7} respectively. For technical reasons, we need a high-order corrector $F_2$ which will be defined in \eqref{C.2.2.2}.

Define the weight function
\begin{align}w(v)=(1+|v|^2)^{\f{\beta}{2}}e^{\varpi|v|^2}\label{w}
\end{align}
with $\beta>3$ and $0<\varpi\leq 1/8$.
The main result in this paper can be stated as follows.

\begin{theorem}\label{thmC.1.1}
Suppose $|\theta_1-\theta_0|\leq \delta_0$ for small $\delta_0$. For sufficiently small $\vep>0$ and any $\a\in (0,\f12)$, there exists a unique solution $F^\vep$ in the form of \eqref{C.1.3.1} to the steady Boltzmann equation \eqref{C.1.0.1} with boundary condition \eqref{C.1.0.2} and total mass condition \eqref{C.2.0.2}. Moreover, there exists constant $p=p(\a)\in (2,\infty)$, such that the remainder term $F_R$ satisfies the following uniform-in-$\vep$ estimate:
	\begin{align}\label{C.1.3.2}
	\left\|\f{wF_R}{\sqrt{\mu}}
	\right\|_{L^\infty}+\vep^{-\f1p}\left\|\left(\f{\P_{\CM}F_R}{\sqrt{\CM}}\right)\right\|_{L^p}
	+\vep^{-1-\f1p}\left\|\f{\nu^{\f12}(I-\P_{\CM})F_R}{\sqrt{\CM}}\right\|_{L^2}\leq C_{\alpha}|\theta_1-\theta_0|.
	\end{align}
	Here the constant $C_\alpha>0$ is uniform in $\vep.$
\end{theorem}
\begin{remark}\label{rmk1}
Esposito et al. in \cite{ELM,ELM2}  studied the hydrodynamic limit of \eqref{C.1.0.1} with \eqref{C.2.0.2}, in the presence of a small external force. They proved that the solution converges to  the steady CNS with {\it no-slip} boundary condition. In this paper, we aim to justify the more accurate CNS approximation by taking account into the slip boundary conditions. Thanks to this choice, we can avoid the higher order expansions used in \cite{ELM,ELM2}.
\end{remark}

The hydrodynamic limit is one of the most fundamental problems in kinetic theory. There are extensive studies on the mathematical description of relations between Boltzmann equation and various of hydrodynamic models. Now we review some of them which are most related to the topic of this paper. For more detailed references, we refer to the book by Cercignani \cite{C} and the survey book by Saint-Raymond \cite{SR}.

Let us first focus on the Euler scaling. The first  mathematical proof of the compressible Euler limit was given by Nishida \cite{N} in the analytic framework. An extension of this result has been made in \cite{UA} for the case when the solution contains initial layers.
By using a truncated Hilbert expansion, Caflisch \cite{Caf} justified the Euler limit for any given smooth Euler solutions; see also \cite{GJJ} for the result in $L^2$-$L^\infty$ framework.  In the same spirit as \cite{Caf}, Lachowicz \cite{L} justified the CNS approximation over the short time interval. Recently, the global-in-time CNS approximation was justified by the second and third authors in a paper  with Zhao \cite{LYZ}
for the case when the data are close to the global equilibrium. This result was extended to case of a general bounded domain in \cite{DL}. On the other hand, the hydrodynamic limit to the compressible Navier-Stokes equations for the steady Boltzmann equation in a slab was studied by Esposito-Lebowitz-Marra \cite{ELM,ELM2}; see also a recent survey \cite{EM}.
We also refer to \cite{HWWY,XZ,Yu} for hydrodynamic limits to some wave patterns. Very recently, the compressible Euler limit in the half-space was studied in \cite{GHW} with  the specular reflection boundary condition.

In diffusive scaling, there are many interesting results on the hydrodynamic limits to the incompressible fluid systems in different settings, {\it cf.} \cite{AEMN,BGL2,BU,EGKM-hy,GS,Guo3,JK,JM,Wu} and the references therein.

The rest of the paper is organized as follows. In Section \ref{sec2} we will present some basic estimates on linear and nonlinear collision terms. In Section \ref{sec3},  the construction of approximate solutions
is given. Precisely, in Section \ref{sec3.1}, we solve the steady Navier-Stokes equation with slip boundary conditions. Some properties of Knudsen layer $\mathfrak{B}_0$ and $\mathfrak{B}_1$ are given in Section \ref{sec3.2}. We  construct the higher order corrector $F_2$ and give some error bounds in Section \ref{sec3.3}. In Section \ref{sec4}, we will study the linearized steady Boltzmann equation. In Section \ref{sec5}, we further construct the remainder $F_R$ and give the proof of Theorem \ref{thmC.1.1}. In Appendix, we summarize some properties of the solution to the Milne problem.

{\it Notations.} Throughout the paper, we denote by $C$ a generic positive constant and by $C_a$ a constant depending on $a$. These constants may vary from line to line. Let $1\leq p\leq \infty$, we denote by $\|\cdot\|_{L^p}$ the $L^p(\Omega\times \mathbb{R}^3)$ norm. We use $|\cdot|_{L^p_v}$ and $|\cdot|_{L^p_x}$ respectively to denote the $L^p(\mathbb{R}^3)$-norm  in the velocity variable and $L^p([0,1])$-norm in the space variable.
For the phase boundary $\g_\pm$ \eqref{C.a.10}, we set measures $\dd\g_{\pm}$ on $\g_\pm$ as
$$\int_{\g_+}f\dd \g_+=\int_{v_1<0}f(0,v)|v_1|\dd v+\int_{v_1>0 }f(1,v)v_1\dd v,
$$
and
$$\int_{\g_-}f\dd \g_-=\int_{v_1>0}f(0,v)v_1\dd v+\int_{v_1<0 }f(1,v)|v_1|\dd v.
$$
For any $1\leq p\leq \infty$, we denote by $|\cdot|_{L^p_{\pm}}$ the $L^p$-norm on $\gamma_\pm$. For $p=2$, we denote by
$\langle \cdot,\cdot\rangle_{\gamma_\pm}$ the inner product on $\gamma_\pm$, that is,
\begin{align}
\langle f,g\rangle_{\gamma_\pm}:=\int_{\g_\pm}f g\dd\g_\pm.\nonumber
\end{align}


\section{Estimates on collision operators \label{sec2}}
Let $\mu(v)=\frac{1}{(2\pi)^{\f32}}e^{-\f{1}{2}|v|^2}$. The linearized collision operator is defined by 
\begin{align}
Lf:=-\f{1}{\sqrt{\mu}}[Q(\mu,\sqrt{\mu}f)+Q(\sqrt{\mu}f,\mu)].\nonumber
\end{align}
As in \cite{Gra},  we have the decomposition $L=\nu-K$, where
$$\nu(v)=\int_{\mathbb{R}^3}\int_{\mathbb{S}^2}B(v-u,\omega)\mu(u)\,\dd\omega\dd u\sim 1+|v|,
$$
and $K=K_1-K_2$ are defined by
\begin{align}
(K_1f)(v)&=\int_{\mathbb{R}^3}\int_{\mathbb{S}^2}B(v-u,\omega)\sqrt{\mu(v)\mu(u)}f(u)\,\dd\omega\dd u,\nonumber\\
(K_2f)(v)&=\int_{\mathbb{R}^3}\int_{\mathbb{S}^2}B(v-u,\omega)\sqrt{\mu(u)\mu(u')}f(v')\,\dd\omega\dd u+\int_{\mathbb{R}^3}\int_{\mathbb{S}^2}B(v-u,\omega)\sqrt{\mu(u)\mu(v')}f(u')\,\dd\omega\dd u.\nonumber
\end{align}
\begin{lemma}[cf. \cite{Gra,Guo2}]
	$K$ is an integral operator given by
	$$
	Kf:=\int_{\mathbb{R}^3}k(v,u)\,\dd u,
	$$
	where
	\begin{align}\label{k1}
	|k(v,u)|\leq C\left\{|v-u| +|v-u|^{-1} \right\}e^{-\f{|v-u|^2}{8}}e^{-\f{||v|^2-|u|^2|^2}{8|v-u|^2}},
	\end{align}
	for any $v, u\in \mathbb{R}^3$ with $v\neq u$.
	Moreover, for the weight function $w(v)$ given by \eqref{w}, it holds that
	\begin{align}\label{k2}
	\int_{\mathbb{R}^3} \left|k(v,u)\right|w^{-1}(u)\dd u \leq (1+|v|)^{-1}w^{-1}(v).
	\end{align}
\end{lemma}
The following lemma gives some estimates on the nonlinear collision operator $Q(f,g)$.
\begin{lemma} [cf. \cite{GPS,Guo2,LYYZ}]
Let $\CM_*$ be any Maxwellian and $w(v)$ be the weight function defined in \eqref{w}. It holds that
	\begin{align}
	\left|\f{\nu^{-1}wQ(f,g)}{\sqrt{\CM_*}}\right|_{L^{\infty}_v}&\leq C\left|\f{wf}{\sqrt{\CM_*}}\right|_{L^{\infty}_v}\left|\f{wg}{\sqrt{\CM_*}}\right|_{L^{\infty}_v},\label{g1}
\end{align}
and
\begin{align}
\left|\f{\nu^{-\f12}Q(f,g)}{\sqrt{\CM_*}}\right|_{L^{2}_v}&\leq C\left|\f{\nu^{\f12}f}{\sqrt{\CM_*}}\right|_{L^2_v}\left|\f{g}{\sqrt{\CM_*}}\right|_{L^{2}_v}+C\left|\f{\nu^{\f12}f}{\sqrt{\CM_*}}\right|_{L^2_v}\left|\f{g}{\sqrt{\CM_*}}\right|_{L^{2}_v}.\label{g2}
	\end{align}
\end{lemma}
\section{Approximate solutions \label{sec3}}
\subsection{Steady Navier-Stokes equations \label{sec3.1}} In this subsection, we construct the solution to the steady Navier-Stokes equations \eqref{C.1.1.3} with slip boundary conditions \eqref{C.1.2.8} and \eqref{C.1.2.9}.  By $\eqref{C.1.1.3}_1$ and boundary condition $u_1(0)=u_1(1)=0$, we have $u_1\equiv0.$ Then by $\eqref{C.1.1.3}_3$ and boundary conditions \eqref{C.1.2.8}, \eqref{C.1.2.9} for  $u_2,u_3$, we have
$u_2,u_3\equiv0.$ Thus, the original problem  \eqref{C.1.1.3},\eqref{C.1.2.8} and \eqref{C.1.2.9}  is reduced to 
\begin{equation}\label{C.2.1.4}
\left\{\begin{aligned}
&\rho \theta\equiv P_0,\\
&\frac{d}{dx}(\sqrt{\theta}\frac{d\theta}{dx})=0,~ x\in(0,1),\\
&\f{\theta(0)-\theta_0}{\theta(0)}=\f{1}{P_0}c_{\beta,2}\vep\frac{d\theta}{dx}(0),~~ \f{\theta(1)-\theta_1}{\theta(1)}=-\f{1}{P_0}c_{\beta,2}\vep\frac{d\theta}{dx}\theta(1),
\end{aligned}\right.
\end{equation}
where $P_0>0$ is a given constant.
\begin{lemma}\label{lmC.2.1.1}
	There exists $\vep_0>0$, such that for any $\vep\in (0,\vep_0)$, there exist a constant $P_0$ and a unique solution $(\rho_{NS},\theta_{NS})$ to \eqref{C.2.1.4} such that
  \begin{align}\label{C.2.1.2-1}\int_0^1\rho_{NS}(x)\dd x=1,
	\end{align}
and
	\begin{align}\label{C.2.1.2}
	\sum_{i=0}^k\left|\frac{d^k}{dx^k}(\rho_{NS}-1,\theta_{NS}-\theta_0)\right|_{L^\infty_x}\leq C_k|\theta_1-\theta_0|,~ \forall k\in \mathbb{N}.
	\end{align}
\end{lemma}
{\bf Proof:}
The general solutions are
	\begin{align}
	\theta(x)=(D_1x+D_2)^{\f23},~\rho(x)=P_0\theta^{-1}(x),\label{C.2.1.2-3}
	\end{align}
	where $D_1$, $D_2$ and $P_0$ are constants to be determined. To satisfy the boundary conditions $\eqref{C.2.1.4}_3$ and total mass condition \eqref{C.2.1.2-1}, we take $(D_1,D_2,P_0)$ as the solution to the following algebraic system:
	\begin{equation}\label{C.2.1.6}\left\{
	\begin{aligned}
	&\CF_1(\vep;D_1,D_2,P_0):=D_2^{\f23}-\theta_0-\f{2c_{\beta,2}}{3P_0}\vep D_1D_2^{\f13}=0,\\
	&\CF_2(\vep;D_1,D_2,P_0):=(D_1+D_2)^{\f23}-\theta_1+\f{2c_{\beta,2}}{3P_0}\vep D_1(D_1+D_2)^{\f13}=0,\\
	&\CF_3(\vep;D_1,D_2,P_0):=\int_{0}^1\rho(x)\dd x-1=\f{3P_0}{D_1}\left[(D_1+D_2)^{\f13}-D_2^{\f13}\right]-1=0.
	\end{aligned}\right.
	\end{equation}
	Notice that when $\vep=0$, \begin{align}
	\left(	D_{1,*},D_{2,*},P_{0,*}\right)\eqdef\left(\theta_1^{\f32}-\theta_0^{\f32},~\theta_0^{\f32},~\f{(\theta_1^{\f32}-\theta_0^{\f32})}{3(\theta^{\f12}_1-\theta^{\f12}_0)}\right)\nonumber
	\end{align}	
		 is the solution of \eqref{C.2.1.6}. By a straightforward calculation, we obtain the Jacobian determinant at $\left(	D_{1,*},D_{2,*},P_{0,*}\right)$ is
	$$\left|\f{\pa(\CF_1,\CF_2,\CF_3)}{\pa(D_1,D_2,P_0)}\right|=\frac{4}{9}\theta_0^{-\f12}\theta_1^{-\f12}P_{0,*}^{-1}+O(1)\vep,
	$$
which does not vanish for any $\vep\in (0,\vep_0)$ with small $\vep_0$. Then by the implicit function theorem, there is a unique solution $[D_1,D_2,P_0]$ of \eqref{C.2.1.6} for any $\vep\in (0,\vep_0)$. The estimate \eqref{C.2.1.2} follows from the explicit formula \eqref{C.2.1.2-3}. The proof of Lemma \ref{lmC.2.1.1} is completed. \qed

\begin{remark}
	The boundary conditions $\eqref{C.2.1.4}_3$ means that there is a temperature gap which is proportional to the normal derivatives of temperature, between fluid layer and the boundary. The proportional coefficient is of the same order as the scale of Knudsen layer.
\end{remark}
\begin{remark}
	Since the pressure $P_0=\rho_{NS}\theta_{NS}$ is a positive constant,  $[\rho_{NS},0,\theta_{NS}]$ is also a solution to steady Euler equations.
\end{remark}
\subsection{Knudsen layers \label{sec3.2}}
In this subsection, we summarize some properties of Knudsen layers $\mathfrak{B}_0$ and $\mathfrak{B}_1$ which are defined in \eqref{C.a.8} and \eqref{C.a.7} respectively.
\begin{lemma}\label{lmC.1.2.1}
	Let $y=x/\vep$ be the stretched variable. Then $\mathfrak{B}_0$ is a solution to the following half-space problem:
	\begin{equation}\left\{
	\begin{aligned}\label{C.1.2.1}
	&v_1\pa_y\mathfrak{B}_0+\rho(0)\sqrt{\frac{2\pi}{\theta_0}}\CL_{\mu_{\theta_0}}B_0=0, y>0,v\in \mathbb{R}^3,\\
	&\mathfrak{B}_0(0,v)|_{v_1>0}=-G_0(0)-\Psi(0),\\
	&\int_{\mathbb{R}^3}v_1\mathfrak{B}_0\dd v\equiv 0,\\
	&\lim_{y\rightarrow +\infty} \mathfrak{B}_0=0,
	\end{aligned}
	\right. 	
	\end{equation}
	where $\Psi_0$ is given by \eqref{C.1.2.5}. Moreover, for any $\varpi\in (0,\f14)$ and $\beta>3$, there exist positive constants $C>0$ and $\sigma_0>0$, such that
	\begin{align}\label{C.1.2.6-2}
		\left|(1+|\cdot|^2)^{\f\beta2}\f{e^{\f{\varpi|\cdot|^2}{\theta_0}}\mathfrak{B}_0(y,\cdot)}{\sqrt{\mu_{\theta_0}}}\right|_{L^\infty_v}\leq Ce^{-\sigma_0 y}|\theta_0-\theta_1|, ~\forall y>0.
	\end{align}
\end{lemma}
{\bf Proof:} From the ansatz in subsection \ref{sp-bcd}, it is direct to check that $\mathfrak{B}_0$ satisfies \eqref{C.1.2.1}. The estimate \eqref{C.1.2.6-2}  follows from the explicit formula \eqref{C.a.8}, \eqref{C.2.1.2} and \eqref{C.a.3} in Lemma \ref{lmC.a.1}. We omit the details for brevity.\qed

Similarly, for $\mathfrak{B}_1$, we have the following lemma.

\begin{lemma}\label{lmC.1.2.2}
	Let $y'=(1-x)/\vep.$ $\mathfrak{B}_1$ satisfies	\begin{equation}\left\{
	\begin{aligned}
	&v_1\pa_{y'}\mathfrak{B}_1+\rho(1)\sqrt{\frac{2\pi}{\theta_1}}\CL_{\mu_{\theta_1}}B_1=0, y'>0,v\in \mathbb{R}^3,\\
	&\mathfrak{B}_1(0,v)|_{v_1>0}=-G_0(1)-\Psi(1),\\
	&\int_{\mathbb{R}^3}v_1\mathfrak{B}_1\dd v\equiv 0,\\
	&\lim_{y'\rightarrow+\infty} \mathfrak{B}_1=0,\nonumber
	\end{aligned}
	\right. 	
	\end{equation}
	where $\Psi(1)$ is defined in \eqref{C.1.2.5-2}. Moreover, $\mathfrak{B}_1$ satisfies the following estimate
	\begin{align}\label{C.1.2.8-2}
		\left|(1+|\cdot|^2)^{\f\beta2}\f{e^{\f{\varpi|\cdot|^2}{\theta_1}}\mathfrak{B}_1(y',\cdot)}{\sqrt{\mu_{\theta_1}}}\right|_{L^\infty_v}\leq Ce^{-\sigma_1 y'}|\theta_0-\theta_1|,~ \forall y>0,
	\end{align}
	for $\varpi\in (0,\f14)$ and $\beta>3$.
\end{lemma}
\subsection{Error terms \label{sec3.3}} Recall  the ansatz \eqref{C.1.3.1}, for a corrector $F_2$ with macroscopic component to  be determined later. Inserting \eqref{C.1.3.1} into \eqref{C.1.0.1} yields the following equation of remainder $F_R$:
\begin{align}\label{C.2.2.3}
v_1\pa_xF_R+\vep^{-1}\CL_{\CM}F_R=-\CL_{as}F_R+\vep^{\alpha}Q(F_R,F_R)+\vep^{-\alpha}A_s,
\end{align}
where
$$\CL_{as}F_R=-Q(G+\MRB_0+\MRB_1+\vep F_2,F_R)-Q(F_R,G+\MRB_0+\MRB_1+\vep F_2).
$$
$A_s$ is related to the error term due to the fact that the solution considered here is an approximation. Moreover, it follows that
\begin{align}A_s=&-\vep^{-1}\left[(\CL_{\CM}-\rho(0)\sqrt{\f{2\pi}{\theta_0}}\CL_{\mu_{\theta_0}})\MRB_0+(\CL_{\CM}-\rho(1)\sqrt{\f{2\pi}{\theta_1}}\CL_{\mu_{\theta_1}})\MRB_1\right]\nonumber\\
&+Q(\MRB_0+\MRB_1,G+\MRB_0+\MRB_1+\vep F_2)+Q(G+\vep F_2, \MRB_0+\MRB_1)-\vep v_1\pa_xF_2.
\label{C.2.2.3-1}
\end{align}
Recall \eqref{F2} for $F_2$. We need to carefully choose the macroscopic component of $F_2$ i.e. $\tilde{F}_2$.
The motivation is twofold. On one hand, notice that there is a possible non-vanishing total mass of the Knudsen layers:
$$m(\vep)=\int_0^1\int_{\mathbb{R}^3}\MRB_0(\frac{x}{\vep},v)+\MRB_1(\frac{1-x}{\vep},v)\dd v\dd x
$$
that is of order $\vep$ by \eqref{C.1.2.6-2} and \eqref{C.1.2.8-2}.
$\tilde{F}_2$ is used to eliminate the extra mass. On the other hand, $\tilde{F}_2$ is chosen so that the residual $A_s$ is purely microscopic. This leads to require
\begin{align}\label{C.2.2.2-2}
\int_0^1\int_{\mathbb{R}^3} F_2(x,v)\dd v\dd x=-\vep^{-1}m(\vep),~~	\P_{\CM}v_1\pa_xF_2=0.
\end{align}
Fortunately, the above requirements can be achieved by choosing \begin{align}\tilde{F}_2=\f{-\vep^{-1}m(\vep) }{\rho_{NS}(x)}\CM+\chi(x)\f{|v|^2-3\theta_{NS}}{2\theta_{NS}}\CM,\nonumber
\end{align}
where \begin{align}\chi(x)=\f{1}{P_0}\left\{\vep^{-1}m(\vep)\theta_{NS}(x)-
\int_{\mathbb{R}^3}v_1^2\CL_{\CM}^{-1}[-(I-P_\CM)v_1\pa_xG+Q(G,G)]\dd v\right\}.\nonumber\end{align}
Then $F_2$ is given by
\begin{align}
F_2=\tilde{F}_2	+L^{-1}_{\CM}[-(I-P_{\CM})v_1\pa_{x}G+Q(G,G)].\label{C.2.2.2}
\end{align}
It is straightforward to check that
	$\int_{\mathbb{R}^3}v_1F_2\dd v=\int_{\mathbb{R}^3}v_1^2F_2\dd v\equiv 0$
and $\int_{\mathbb{R}^3}v_1|v|^2F_2\dd v\equiv0$ since $G$ is odd in $v_1$. Therefore, \eqref{C.2.2.2-2} holds.

The boundary condition of $F_R$ is given by
\begin{equation}
\begin{aligned}\label{C.2.2.4}
F_R|_{\g_-}=\mu_{w}\int_{\{n(x)\cdot v>0\}}F_R\{n(x)\cdot v\}\dd v+\vep^{1-\alpha}r.
\end{aligned}
\end{equation}
Here
$$\begin{aligned}&r(0,v)=\vep^{-\alpha}\left(\MRB_{1}|_{x=0}-\mu_{\theta_0}\int_{\{v_1<0\}}\MRB_1|_{x=0}|v_1|\dd v\right)\\
&\quad+M_R(0)+ G_R(0)+ F_2(0)-\mu_{\theta_0}\int_{\{v_1<0\}}[M_R(0)+ G_R(0)+F_2(0)]v_1\dd v,
\end{aligned}
$$
and
$$\begin{aligned}&r(1,v)=\vep^{-\alpha}\left(\MRB_{0}|_{x=1}-\mu_{\theta_1}\int_{\{v_1>0\}}\MRB_0|_{x=1}|v_1|\dd v\right)\\
&\quad+ M_R(1)+ G_R(1)+ F_2(1)-\mu_{\theta_1}\int_{\{v_1>0\}}[ M_R(1)+G_R(1)+ F_2(1)]v_1\dd v,
\end{aligned}
$$
where $M_R$ and $G_R$ are defined in \eqref{C.1.2.2} and \eqref{C.1.2.2-1} respectively.

We conclude this section by summarizing some estimates on approximate solutions and errors for later use.
\begin{lemma}\label{lem2.1}	Let $w$ be the weight function defined in \eqref{w}. We have
\begin{align}
&|G(x,v)|\leq C (1+|v|)^4\CM|\pa_x\theta_{NS}(x)|,~\left\|\f{wF_2}{\sqrt{\mu}}\right\|_{L^\infty}+\left\|\f{w\pa_xF_2}{\sqrt{\mu}}\right\|_{L^\infty}\leq C,\label{6.2}\\
&\left|\frac{w\mathfrak{B}_0}{\sqrt{\mu}}\right|_{L^\infty_v}\leq Ce^{-\sigma_1\frac{x}{\vep}}|\pa_x\theta_{NS}|_{L^\infty_x},~ \left|\frac{w\mathfrak{B}_1}{\sqrt{\mu}}\right|_{L^\infty_v}\leq Ce^{-\sigma_1\frac{1-x}{\vep}}|\pa_{x}\theta_{NS}|_{L^\infty_x},\label{6.3}\\
&\left\|\frac{\nu^{-1}wA_s}{\sqrt{\mu}}\right\|_{L^\infty}\leq C|\pa_x\theta_{NS}|_{L^\infty_x},~\left\|\frac{\nu^{-\f12}A_s}{\sqrt{\CM}}\right\|_{L^2}\leq C\sqrt{\vep}|\pa_x\theta_{NS}|_{L^\infty_x},\label{6.4}\\
&\left|\frac{wr}{\sqrt{\mu}}\right|_{L^\infty_-}+\left|\frac{r}{\sqrt{\CM}}\right|_{L^2_-}\leq C|\pa_x\theta_{NS}|_{L^\infty_x}.\label{6.5}
\end{align}
\end{lemma}
{\bf Proof:} The estimates \eqref{6.2}-\eqref{6.5} are straightforward by using the explicit formula and bounds \eqref{C.1.2.6-2}, \eqref{C.1.2.8-2}. Recall $A_s$ defined in \eqref{C.2.2.3-1}. Consider the highest singular term  $\CQ=\vep^{-1}(\CL_{\CM}-\rho(0)\sqrt{\frac{2\pi}{\theta_0}}\CL_{\mu_{\theta_0}})\mathfrak{B}_0.$ By mean value theorem, one has $$\begin{aligned}
|\CM-\rho(0)\sqrt{\frac{2\pi}{\theta_0}}\mu_{\theta_0}|&\leq 	|\CM(x)-\CM(0)|
+\left|\CM(0)-\rho(0)\sqrt{\frac{2\pi}{\theta_0}}\mu_{\theta_0}\right|\\
&\leq C\left(|\pa_x\theta_{NS}|_{L^\infty_x}|x|+|\theta_{NS}(0)-\theta_0|\right)(1+|v|)^2\mu_{\theta_0},\\
&\leq C|\pa_x\theta_{NS}|_{L^\infty_x}\left(|x|+\vep\right)(1+|v|)^2\mu_{\theta_0},\\
\end{aligned}$$
where we have used the boundary condition \eqref{C.2.1.4} for $\theta_{NS}(0)$ in the last inequality. Then by \eqref{g2} and \eqref{C.1.2.6-2}, it holds that
\begin{equation}
\begin{aligned}
\left|\f{\nu^{-\f12}\CQ}{\sqrt{\CM}}\right|_{L^2_{v}}
&\leq C\vep^{-1}
\left|\f{\nu^{\f12}\left(\CM-\rho(0)\sqrt{\frac{2\pi}{\theta_0}}\mu_{\theta_0}\right)}{\sqrt{\CM}}\right|_{L^2_v}\left|\f{\mathfrak{B}_0}{\sqrt{\CM}}\right|_{L^2_v}\nonumber\\
&\quad+C\vep^{-1}\left|\f{\CM-\rho(0)\sqrt{\frac{2\pi}{\theta_0}}\mu_{\theta_0}}{\sqrt{\CM}}\right|_{L^2_v}\left|\f{\nu^{\f12}\mathfrak{B}_0}{\sqrt{\CM}}\right|_{L^2_v}\\
&\leq C|\pa_x\theta_{NS}|_{L^\infty_x}\left(\left|\f{x}{\vep}\right|+1\right)e^{-\sigma_0\left|\f{x}{\vep}\right|}.
\end{aligned}
\end{equation}
This implies that
$$\left\|\f{\nu^{-1/2}\CQ}{\sqrt{\CM}}\right\|_{L^2}\leq C\sqrt{\vep}|\pa_x\theta_{NS}|_{L^\infty_x}.
$$
Other terms can be estimated similarly and we omit the details for brevity. The proof of Lemma \ref{lem2.1} is completed. \qed

\section{Linear problem \label{sec4}}
In this section, we will study the following linear stationary problem:
\begin{equation}\label{C.2.3.1}
\left\{
\begin{aligned}
&\vep v_1\pa_xF_R+\CL_{\CM}F_R=g,\quad x\in (0,1), v\in \mathbb{R}^3,\\
&F_R|_{\g_-}=\P_\g F_R+r,
\end{aligned}\right.
\end{equation}
where 
\begin{align}\nonumber
\left(\P_\g f\right)(x,v)=\frac{1}{2\pi\theta_{NS}^2(x)}e^{-\frac{|v|^2}{2\theta_{NS}(x)}}\int_{n(x)\cdot u>0}f(x,u)\{n(x)\cdot u\}\dd u,
\end{align}
and both $g$ and $r$ are inhomogeneous source terms. The following is the main result in this section.
\begin{proposition}\label{propC.2.1} Assume that
	\begin{align}\label{C.2.3.2}
	\int_{\mathbb{R}^3}g(x,v)\dd v=0,~ \forall x\in (0,1),
	\end{align}
	and
	\begin{align}\label{C.2.3.3}
	\int_{\{v_1>0\}}r(0,v)v_1\dd v=\int_{\{v_1<0\}}r(1,v)v_1\dd v=0.
	\end{align}
	There exists a positive constant $\delta_0>0$, such that if $|\theta_1-\theta_0|\leq \delta_0$, then for any sufficiently small $\vep>0$ and for any $p\in (2,\infty)$, the linear problem \eqref{C.2.3.1} admits a unique solution $F_R$ satisfying
	\begin{align}\label{C.2.3.4-1}
	\int_{0}^1\int_{\mathbb{R}^3}F_R(x,v)\dd x\dd v=0,
	\end{align}
	and
	\begin{align}\label{C.2.3.4}
	\left\|\f{wF_R}{\sqrt{\mu}}\right\|_{L^\infty}&+\vep^{-\f1p}\left\|\f{\P_\CM F_R}{\sqrt{\CM}}\right\|_{L^p}+\vep^{-1-\f1p}\left\|\f{\nu^{\f12}(I-\P_\CM)F_R}{\sqrt{\CM}}\right\|_{L^2}\nonumber\\
	&\leq  C_p\vep^{-2-\f1p}\left\|\f{\P_{\CM}g}{\sqrt{\CM}}\right\|_{L^2}+C_p\vep^{-1-\f1p}\left\|\f{\nu^{-\f12}(I-\P_\CM)g}{\sqrt{\CM}}\right\|_{L^2}\nonumber\\
	&\qquad+C_p\vep^{-\f12-\f1p}\left|\f{r}{\sqrt{\CM}}\right|_{L^2_{-}}+C_p\left\|\f{\nu^{-1}wg}{\sqrt{\mu}}\right\|_{L^\infty}+C_p\left|\f{wr}{\sqrt{\mu}}\right|_{L^\infty_-}.	\end{align}
 Here, the positive constant $C_p>0$ does not depend on $\vep.$
\end{proposition}

\subsection{$L^2$-estimate}
\begin{lemma}\label{lem3.1}
	Let $F_R$ be a solution to the linearized problem \eqref{C.2.3.1}. Assume that $g$ and $r$ satisfy \eqref{C.2.3.2} and \eqref{C.2.3.3} respectively. There exists $\tau_0>0$, such that for sufficiently small $\vep\ll1$ and for any
	$\eta\in (0,1)$, it holds
	\begin{align}\label{leme1}
&\left\|\frac{\nu^{\f12}(I-\P_{\CM})F_R}{\sqrt{\CM}}\right\|_{L^2}+\vep^{\f12}\left|\f{(I-\P_\g)F_R}{\sqrt{\CM}}\right|_{L^2_+}\nonumber\\
&\quad\leq \vep\left(\eta+|\pa_x\theta_{NS}|_{L^\infty_x}\right)\left\|\f{\P_{\CM}F_R}{\sqrt{\CM}}\right\|_{L^2}+C_\eta e^{-\f{\tau_0}{\vep}}\left\|\f{wF_R}{\sqrt{\mu}}\right\|_{L^\infty}\nonumber\\
&\qquad +C_{\eta}\left(\vep^{-1}\left\|\f{\P_\CM g}{\sqrt{\CM}}\right\|_{L^2}+\left\|\f{\nu^{-\f12}(I-\P_\CM)g}{\sqrt{\CM}}\right\|_{L^2}+\vep^{\f12}\left|\f{r}{\sqrt{\CM}}\right|_{L^2_-}
\right).
\end{align}
Here, the constant $C_\eta$ does not depend on $\vep$.
\end{lemma}
{\bf Proof:} For simplicity, we denote $$\CM_\omega=\frac{1}{2\pi\theta_{NS}^2}e^{-\frac{|v|^2}{2\theta_{NS}}}.$$
By taking inner product of \eqref{C.2.3.1} with $F_R/\CM_\omega$ over $(0,1)\times \mathbb{R}^3$ and then integrating by parts, we have
	\begin{align}\label{C.2.3.6}
	&\f\vep{2}\int_{\g_+}\frac{F_R^2}{\CM_\omega}\dd \gamma_+-\f\vep{2}\int_{\g_-}\frac{F_R^2}{\CM_\omega}\dd \gamma_-+\int_0^1\int_{\mathbb{R}^3}\frac{F_R\CL_\CM F_R}{\CM_\omega}\dd x\dd v\nonumber\\
	&\qquad+\int_0^1\int_{\mathbb{R}^3}\frac{v_1\pa_x\CM_\omega}{2\CM^2_\omega}|F_R|^2\dd x\dd v=\int_0^1\int_{\mathbb{R}^3}\frac{g F_R}{\CM_\omega}\dd x\dd v.
	\end{align}
By \eqref{C.2.3.3}, we have
\begin{align}
\int_{\ga_-}\frac{r\cdot \P_\g F_R}{\CM_{\omega}}\dd \g_-=&\left(\int_{\{v_1>0\}}r(0,v)v_1\dd v\right)\times \left(\int_{\{v_1<0\}}F_R(0,v)|v_1|\dd v\right)\nonumber\\
&+\left(\int_{\{v_1<0\}}r(1,v)|v_1|\dd v\right)\times \left(\int_{\{v_1<0\}}F_R(1,v)v_1\dd v\right)=0.\nonumber
\end{align}
Thus, by boundary condition $\eqref{C.2.3.1}_2$, we can obtain
	$$\begin{aligned}\f12\int_{\g_-}\frac{F_R^2}{\CM_\omega}\dd \gamma_-&=\f12\int_{\g_-}\frac{|\P_{\g}F_R|^2}{\CM_\omega}\dd \gamma_-+\int_{\g_-}\frac{r\cdot\P_{\g}F_R}{\CM_\omega}\dd \gamma_-+\f12\int_{\g_-}\frac{|r|^2}{\CM_\omega}\dd \gamma_-\\
	&=\f12\int_{\g_-}\frac{|\P_{\g}F_R|^2}{\CM_\omega}\dd \gamma_++\f12\int_{\g_-}\frac{|r|^2}{\CM_\omega}\dd \gamma_-\\
	&=\f12\int_{\g_+}\frac{|\P_{\g}F_R|^2}{\CM_\omega}\dd \gamma_++\f12\int_{\g_-}\frac{|r|^2}{\CM_\omega}\dd \gamma_-.
	\end{aligned}$$
	Hence,  it holds that
	\begin{align}\label{C.2.3.7}\f\vep{2}\int_{\g_+}\frac{F_R^2}{\CM_\omega}\dd \gamma_+-\f\vep{2}\int_{\g_-}\frac{F_R^2}{\CM_\omega}\dd \gamma_-&=\f\vep{2}\int_{\g_+}\frac{|F_R|^2-|\P_\g F_R|^2}{\CM_\omega}\dd \gamma_+-\f\vep{2}\int_{\g_-}\frac{|r|^2}{\CM_\omega}\dd \gamma_-.\nonumber\\
	&=\f\vep{2}\int_{\g_+}\frac{|(I-\P_\g)F_R|^2}{\CM_\omega}\dd \gamma_+-\f\vep{2}\int_{\g_-}\frac{|r|^2}{\CM_\omega}\dd \gamma_-.\end{align}
	Note that $\CM_\omega=\frac{\CM}{\rho_{NS}\theta_{NS}^{\f12}}$. Then
	by coercivity estimate \eqref{C.a.12}, we have
	\begin{align}\label{C.2.3.8}
	\int_0^1\int_{\mathbb{R}^3}\frac{F_R\CL_\CM F_R}{\CM_\omega}\dd x\dd v&\geq \inf\left\{\rho_{NS}\theta_{NS}^{\f12}\right\}\int_0^1\int_{\mathbb{R}^3}\frac{F_R\CL_\CM F_R}{\CM}\dd x\dd v\nonumber\\
	&\geq c_1\left\|\f{\nu^{\f12}(I-\P_{\CM})F_R}{\sqrt{\CM}}\right\|_{L^2}^2,
	\end{align}
	for some positive constant $c_1>0$ independent of $\vep$. By Cauchy-Schwarz and Young's inequalities, it holds for any $\eta\in (0,1)$ and any $\ka\in (0,1)$
	\begin{align}\label{C.2.3.9}
	\left|\int_0^1\int_{\mathbb{R}^3}\frac{g F_R}{\CM_\omega}\dd x\dd v \right|\nonumber\leq & \ka\left\|\f{\nu^{\f12}(I-\P_{\CM})F_R}{\sqrt{\CM}}\right\|_{L^2}^2+C_\ka\left\|\f{\nu^{-\f12}(I-\P_{\CM})g}{\sqrt{\CM}}\right\|_{L^2}^2\\
	&+\eta\vep^2\left\|\frac{\P_{\CM}F_R}{\sqrt{\CM}}\right\|_{L^2}^2+C_\eta\vep^{-2}\left\|\f{\P_{\CM}g}{\sqrt{\CM}}\right\|_{L^2}^2.
	\end{align}
 For the last term on the left hand side of \eqref{C.2.3.6}, we divide  it into the following three parts:
	\begin{align}\label{C.2.3.10}
	&\vep\int_0^1\int_{\mathbb{R}^3}\f{v_1\pa_x\CM_\omega}{2\CM_\omega^2}| F_R|^2\dd x\dd v\nonumber\\
	&=\vep\int_0^1\int_{\mathbb{R}^3}\f{v_1\pa_x\CM_\omega}{2\CM_\omega^2}|\P_\CM F_R|^2\dd x\dd v+2\vep\int_0^1\int_{\mathbb{R}^3}\f{v_1\pa_x\CM_\omega}{2\CM_\omega}\P_{\CM}F_R(I-\P_{\CM})f_R\dd x\dd v\nonumber\\
	&\qquad+\vep\int_0^1\int_{\mathbb{R}^3}\f{v_1\pa_x\CM_\omega}{2\CM_\omega^2}|(I-\P_{\CM})F_R|^2\dd x\dd v\nonumber\\
	\quad:&=I_{1}+I_2+I_3.
	\end{align}
	For $I_1$, by integrating \eqref{C.2.3.1} over $\mathbb{R}^3$ and using \eqref{C.2.3.2}, we have
	\begin{align}
	\f{\dd}{\dd x}\int_{\mathbb{R}^3}v_1F_R(x,v)\dd v=\vep^{-1}\int_{\mathbb{R}^3}g(x,v)\dd v=0.\nonumber
	\end{align}
Then by \eqref{C.2.3.3}, we deduce that
\begin{align}
\int_{\mathbb{R}^3}v_1F_R(x,v)\dd v&\equiv \int_{\mathbb{R}^3}v_1F_R(0,v)\dd v=
\int_{\{v_1>0\}}r(0,v)v_1\dd v=0,\nonumber
\end{align}
	which implies $\mathbb{P}_\CM F_R$ is even in $v_1.$ Then it holds 
	\begin{align}\label{C.2.3.10-1}
	I_1=0.
	\end{align}
	Now we fix $\ka\in (0,1)$.
	For $I_2$, by Cauchy-Schwarz and Young's inequalities, we deduce that
	\begin{align}
	|I_2|\leq \ka\left\|\f{\nu^{\f12}(I-\P_{\CM})f_R}{\sqrt{\CM}}\right\|_{L^2}^2+C_\ka\vep^{2}|\pa_x\theta_{NS}|_{L^\infty_x}^2
\left\|\frac{\P_{\CM}F_R}{\sqrt{\CM}}\right\|_{L^2}^2.\label{C.2.3.10-2}
	\end{align}
	For $I_3$, as in \cite{GJ}, we write it as
	$$
	\begin{aligned}I_3\leq& C\vep|\pa_x\theta_{NS}|_{L^\infty_x} \int_0^1\int_{ \mathbb{R}^3}\frac{(1+|v|)^3|(I-\P_\CM)F_R|^2}{\CM}\dd x\dd v\\
	=& C\vep|\pa_x\theta_{NS}|_{L^\infty_x} \int_0^1
	\int_{\{|v|\leq \ka\vep^{-\f12}\}}\f{(1+|v|)^3|(I-\P_\CM)F_R|^2}{\CM}\dd v\dd x
	\\
	&+C\vep|\pa_x\theta_{NS}|_{L^\infty_x}\int_0^1\int_{\{|v|>\ka\vep^{-\f12}\}}\f{(1+|v|)^3|(I-\P_\CM)F_R|^2}{\CM}\dd v\dd x.
\end{aligned}	
$$
For $|v|\leq \ka \vep^{-\f12},$	we have
\begin{align}
(1+|v|)^2\leq C(1+|v|^2)\leq C\vep^{-1}(\vep+\k^2).\nonumber
\end{align}
Then it holds that
\begin{align}
\int_0^1
\int_{\{|v|\leq \ka\vep^{-\f12}\}}\f{(1+|v|)^3|(I-\P_\CM)F_R|^2}{\CM}\dd v\dd x\leq C\vep^{-1}(\vep+\ka^2)\left\|\f{\nu^{\f12}(I-\P_\CM)F_R}{\sqrt{\CM}}\right\|_{L^2}^2.\nonumber
\end{align}
Recall the weight function $w$ defined in \eqref{w}.
For $|v|>\ka \vep^{-\f12}$,  by noticing that \begin{align}
\theta_{NS}(x)\geq \theta_{NS}(0)=1+O(1)\vep,\nonumber
\end{align}
we have
\begin{align}
w^{-1}(1+|v|)^4\sqrt{\f{\mu}{\CM}}&\leq e^{-(\varpi-C\vep)|v|^2}(1+|v|)^{-(\beta-4)}\nonumber\\
&\leq C e^{-\f{\varpi}{2}|v|^2}\leq Ce^{-\frac{\ka^2\varpi}{2\vep}}.\nonumber
\end{align}
Then we can obtain
\begin{align}
\int_0^1\int_{\{|v|>\ka\vep^{-\f12}\}}\f{(1+|v|)^3|(I-\P_\CM)F_R|^2}{\CM}\dd v\dd x&\leq Ce^{-\f{\ka^2\varpi}{\vep}}\left\|\f{wF_R}{\sqrt{\mu}}\right\|_{L^\infty}^2\int_{\mathbb{R}^3}(1+|v|)^{-5}\dd v\nonumber\\
&\leq Ce^{-\f{\ka^2\varpi}{\vep}}\left\|\f{wF_R}{\sqrt{\mu}}\right\|_{L^\infty}^2.\nonumber
\end{align}
Combining these two estimates, we get
\begin{align}\label{C.2.3.10-3}
I_3\leq C(\vep+\ka^2)\left\|\f{\nu^{\f12}(I-\P_\CM)F_R}{\sqrt{\CM}}\right\|_{L^2}^2+Ce^{-\f{\ka^2\varpi}{\vep}}\left\|\f{wF_R}{\sqrt{\mu}}\right\|_{L^\infty}^2.
\end{align}
	Putting estimates \eqref{C.2.3.6}-\eqref{C.2.3.10-3} together gives
	\begin{align}\label{C.2.3.11}
	&\left\|\frac{\nu^{\f12}(I-\P_{\CM})F_R}{\sqrt{\CM}}\right\|_{L^2}^2+\vep\left|\f{(I-\P_\g)F_R}{\sqrt{\CM}}\right|_{L^2_+}^2\nonumber\\
	&\quad\leq
	C(\vep+\ka)\left\|\frac{\nu^{\f12}(I-\P_{\CM})F_R}{\sqrt{\CM}}\right\|_{L^2}^2+
	 C\vep^2\left(\eta+|\pa_x\theta_{NS}|_{L^\infty_x}^2\right)\left\|\f{\P_{\CM}F_R}{\sqrt{\CM}}\right\|_{L^2}^2+C e^{-\f{\ka^2\varpi}{\vep}}\left\|\f{wF_R}{\sqrt{\mu}}\right\|_{L^\infty}^2\nonumber\\
	&\qquad +C_{\eta}\vep^{-2}\left\|\f{\P_\CM g}{\sqrt{\CM}}\right\|_{L^2}^2+C_\ka\left\|\f{\nu^{-\f12}(I-\P_\CM)g}{\sqrt{\CM}}\right\|_{L^2}^2+	C	\vep\left|\f{r}{\sqrt{\CM}}\right|_{L^2_-}^2.
	\end{align}
Here, $\ka\in (0,1)$ and $\eta\in (0,1)$ are two arbitrary constants. By taking $\ka$ suitably small in \eqref{C.2.3.11}, we obtain \eqref{leme1}. The proof of Lemma \ref{lem3.1} is completed.\qed
 \subsection{$L^p$-estimates on $\P_{\CM}F_R$.}

\begin{lemma}\label{lem3.2}
for any $p\in [2,\infty)$, there exists a positive constant $C$, such that
	\begin{align}\label{C.2.3.23}
	\vep\left\|\f{\P_{\CM}F_R}{\sqrt{\CM}}\right\|_{L^p}\leq&C\vep|\pa_x\theta_{NS}|_{L^\infty_x}\left\|\f{\P_{\CM}F_R}{\sqrt{\CM}}\right\|_{L^p}+ C\vep^{\f{p}{p-2}}\left\|\f{wF_R}{\sqrt{\mu}}\right\|_{L^\infty}+C\vep\left|\f{(I-\P_\g)F_R}{\sqrt{\CM}}\right|_{L^2_+}\nonumber\\
	&+C\left\|\f{(I-\P_\CM)F_R}{\sqrt{\CM}}\right\|_{L^2}+C\left\|\f{\nu^{-\f12}g}{\sqrt{\CM}}\right\|_{L^2}+C\vep\left|
	\f{r}{\sqrt{\CM}}\right|_{L^2_-}.
	\end{align}	
	\end{lemma}
{\bf Proof:}	
The proof is based on the dual argument developed in \cite{EGKM,EGKM-hy}. We need to slightly modify their method here since the reference Maxwellian $\CM$ depends on $x$. We decompose $F_R$ as
\begin{align}
F_R&=\P_\CM F_R+(I-\P_\CM)F_R\nonumber\\
&=\frac{a(x)}{\rho_{NS}}\CM+\frac{b(x)\cdot v}{\rho_{NS}\theta_{NS}}\CM+\frac{c(x)}{2\rho_{NS}\theta_{NS}}\left(\frac{|v|^2}{\theta_{NS}}-3\right)\CM+(I-\P_\CM)F_R.\nonumber
\end{align}
Thus, it suffices to obtain $L^p$-estimate of $(a,b,c)$. For simplicity, we only estimate $a$ because the estimates for $b$ and $c$ are similar. Let $\phi(x,v)$ be any smooth test function. Then taking inner product of $\eqref{C.2.3.1}_1$ with $\phi$, we obtain
	\begin{align}\label{C.2.3.12}
	\vep\left[\int_{\mathbb{R}^3}v_1F_R\phi(1,v)-\int_{\mathbb{R}^3}v_1F_R\phi(0,v)\dd v\right]-\vep\int_0^1\int_{ \mathbb{R}^3}v_1F_R\pa_x\phi(x,v)\dd v\dd x\nonumber\\
	\quad+\int_0^1\int_{ \mathbb{R}^3}\left(\CL_{\CM}F_R-g\right)(x,v)\phi \dd v\dd x=0.
	\end{align}
	We take the test function $\phi_a=\pa_x\psi_a\f{v_1(|v|^2-10\theta_{NS})}{\theta_{NS}^2},$ where $\psi_a$ solves
	\begin{equation}\label{C.2.3.17}
	\left\{
	\begin{aligned}
	-\pa_{xx}\psi_a&=a^{p-1}-\int_0^1a^{p-1}(\tau)\dd \tau,\\
	\pa_x\psi_a(0)&=\pa_x\psi_a(1)=0.
	\end{aligned}\right.
	\end{equation}
One can check that
	\begin{align}
	\psi_a=-\int_0^x\left[a^{p-1}(s)-\int_0^1a^{p-1}(\tau)\dd\tau \right](x-s)\dd s\nonumber
	\end{align}
is a solution to \eqref{C.2.3.17} and it satisfies the following estimate
	\begin{align}\label{C.2.3.18-1}
	|\psi_a|_{C^1([0,1])}+|\pa_{xx}\psi_a|_{L^{\f{p}{p-1}}_x}\leq C|a|_{L^p_x}^{p-1}.
	\end{align}
Now we insert $\phi_\a$ into \eqref{C.2.3.12}.	
	Noting that \begin{align}
	\int_{\mathbb{R}^3}v_1^2(|v|^2-10\theta_{NS})(|v|^2-3\theta_{NS})\CM\dd v=	\int_{\mathbb{R}^3}v_1^3(|v|^2-10\theta_{NS})\CM\dd v=0,\nonumber
	\end{align}
	we have
	\begin{equation}\label{C.2.3.18}
	\begin{aligned}
	&\vep\int_0^1\int_{ \mathbb{R}^3}v_1\P_{\CM}F_R\pa_x\phi_a\dd v\dd x\\
	&=\vep\int_0^1\pa_x^2\psi_\a\P_\CM F_R \frac{|v|^2-10\theta_{NS}}{\theta_{NS}^2}v_1^2\dd v\dd x+\vep\int_0^1\int_{ \mathbb{R}^3}v_1^2\P_{\CM}F_R\pa_x\psi_a\pa_x\left(\f{|v|^2-10\theta_{NS}}{\theta^2_{NS}}\right)\dd v\dd x\\
	&=5\vep \int_{0}^1a\left(a^{p-1}-\int_0^1a^{p-1}(\tau)\dd\tau\right)\dd x+\vep\int_0^1\int_{ \mathbb{R}^3}v_1^2\P_{\CM}F_R\pa_x\psi_a\pa_x\left(\f{|v|^2-10\theta_{NS}}{\theta^2_{NS}}\right)\dd v\dd x\\
	&\geq 5\vep |a|_{L^p_x}^p-C\vep|\pa_x\theta_{NS}|_{L^\infty_x}|[a,b,c]|_{L^p_x}|\pa_x\psi_a|_{L^\infty_x}\\
	&\geq 5\vep |a|_{L^p_x}^p-C\vep|\pa_x\theta_{NS}|_{L^\infty_x}|[a,b,c]|_{L^p_x}^p.
	\end{aligned}
	\end{equation}
	Here,  we have used  \eqref{C.2.3.18-1} and \eqref{C.2.3.4-1} so that
	$$\int_{0}^1a\dd x\cdot\int_0^1a^{p-1}(\tau)\dd\tau=\int_{0}^1\int_{\mathbb{R}^3}F_{R}\dd x\dd v\cdot\int_0^1a^{p-1}(\tau)\dd\tau=0.
	$$
	Thanks to the Neumann boundary condition $\eqref{C.2.3.17}_2$,
	the boundary contribution in \eqref{C.2.3.12} vanishes, that is,
	\begin{align}\label{C.2.3.18-2}
	\int_{\mathbb{R}^3}v_1f_R\phi_a(1,v)\dd v=\int_{\mathbb{R}^3}v_1f_R\phi_a(0,v)\dd v=0.
\end{align}
 By H\"{o}lder's inequality, the rest terms in \eqref{C.2.3.12} are bounded as follows:
\begin{align}\label{C.2.3.18-3}
	\vep\left|\int_0^1\int_{ \mathbb{R}^3}v_1(I-\P_\CM)F_R\pa_x\phi_a\dd v\dd x\right|&\leq C\vep \left\|\f{(I-\P_\CM)F_R}{\sqrt{\CM}}\right\|_{L^p}(|\pa_{xx}\psi_a|_{L^{\f{p}{p-1}}_x}+|\pa_x\psi_\a|_{L^\infty_x})\nonumber\\
	&\leq C\vep \left\|\f{(I-\P_\CM)F_R}{\sqrt{\CM}}\right\|_{L^p}|a|_{L^p_x}^{p-1},
	\end{align}
	and
	\begin{align}\label{C.2.3.18-4}
	\left|\int_0^1\int_{ \mathbb{R}^3}\left[\CL_\CM F_R-g\right]\phi_a \dd v\dd x\right|
	&\leq C\left[\left\|\f{(I-\P_{\CM})F_R}{\sqrt{\CM}}\right\|_{L^2}+\left\|\f{\nu^{-\f12}g}{\sqrt{\CM}}\right\|_{L^2}\right]|\pa_x\psi_\a|_{L^\infty_x}\nonumber\\
	&\leq C\left[\left\|\f{(I-\P_{\CM})F_R}{\sqrt{\CM}}\right\|_{L^2}+\left\|\f{\nu^{-\f12}g}{\sqrt{\CM}}\right\|_{L^2}\right]|a|_{L^p_x}^{p-1}.
	\end{align}
Substituting \eqref{C.2.3.18}-\eqref{C.2.3.18-4} into \eqref{C.2.3.12}, we deduce that
	\begin{align}
	\vep|a|_{L^p_x}\leq& C\vep|\pa_x\theta_{NS}|_{L^\infty_x}|[a,b,c]|_{L^p_x}\nonumber\\
	&+C\vep\left\|\f{(I-\P_{\CM})F_R}{\sqrt{\CM}}\right\|_{L^p}+C\left\|\f{(I-\P_{\CM})F_R}{\sqrt{\CM}}\right\|_{L^2}+C\left\|\f{\nu^{-\f12}g}{\sqrt{\CM}}\right\|_{L^2}.\nonumber
	\end{align}
The $L^p$-estimates for $b$ and $c$ can be obtained in the same way. To summarize, we have
\begin{align}\label{C.2.3.19}
\vep\left\|\f{\P_\CM F_R}{\sqrt{\CM}}\right\|_{L^p}\leq & C\vep |\pa_x\theta_{NS}|_{L^\infty_x}\left\|\f{\P_\CM F_R}{\sqrt{\CM}}\right\|_{L^p}+C\vep\left\|\f{(I-\P_{\CM})F_R}{\sqrt{\CM}}\right\|_{L^p}\nonumber\\
&+C\left\|\f{(I-\P_{\CM})F_R}{\sqrt{\CM}}\right\|_{L^2}+C\left\|\f{\nu^{-\f12}g}{\sqrt{\CM}}\right\|_{L^2}.
\end{align}
Combining \eqref{C.2.3.19} with an interpolation inequality
\begin{align}\vep\left\|\f{(I-\P_\CM)F_R}{\sqrt{\CM}}\right\|_{L^p}\leq C\left\|\f{(I-\P_\CM)F_R}{\sqrt{\CM}}\right\|_{L^2}+C\vep^{\f{p}{p-2}}\left\|\f{wF_R}{\sqrt{\mu}}\right\|_{L^\infty}\nonumber
\end{align}
yields \eqref{C.2.3.23}.  Therefore, the proof of Lemma \ref{lem3.2} is completed. \qed

\subsection{Weighted $L^\infty$ estimate}
Recall the normalized  global Maxwellian $\mu(v)=\frac{1}{(2\pi)^{
\f32}}e^{-\f{|v|^2}{2}}$ and associated linearized collision operator $L=\nu-K$. Define $h=\f{wF_R}{\sqrt{\mu}}$. Then the equation of $h$ reads
	\begin{equation}\label{C.2.3.24}
	\left\{
	\begin{aligned}
	&v_1\pa_xh+\vep^{-1}\nu(v)h-\vep^{-1}K_{w}h=\vep^{-1}J,\quad x\in (0,1), v\in \mathbb{R}^3,\\
	& h|_{\g_-}=\f{1}{\tilde{w}(v)}\int_{\{n(x)\cdot v>0\}}h(x,u)\tilde{w}(u)\dd \sigma+q,
	\end{aligned}\right.
	\end{equation}
where we have used the notations:
$$\tilde{w}(v)=\f{1}{w(v)\sqrt{2\pi\mu(v)}},~\dd \sigma=\sqrt{2\pi}\mu(v)\{n(x)\cdot v\}\dd v,~K_wh=wK\left(\f{h}{w}\right).
	$$
$J$ and $q$ are inhomogeneous source which are given by
\begin{align}
J&=\frac{w}{\sqrt{\mu}}\left[Q(\CM-\mu,\f{\sqrt{\mu}h}{w})+Q(\f{\sqrt{\mu}h}{w},\CM-\mu)\right]+\frac{wg}{\sqrt{\mu}},\label{C.3.1}\\
q&=\frac{w}{\sqrt{\mu}}\left(\f{1}{2\pi\theta_{NS}^2}e^{-\f{|v|^2}{2\theta_{NS}}}-\frac{1}{2\pi}e^{-\frac{|v|^2}{2\pi}}\right)\left(\int_{\{n(x)\cdot v>0\}}h(x,u)\tilde{w}(u)\dd \sigma\right)+\frac{wr}{\sqrt{\mu}}.\label{C.3.2}
\end{align}
Fix $t>0$ as a parameter. Given any $(x,v)\in (0,1)\times \mathbb{R}^3$, let  $[X(s),V(s)]$ be the backward bi-characteristics, which is determined by
\begin{equation}
\left\{
\begin{aligned}
&\frac{\dd X(s)}{\dd s}=V^1(s),~\frac{\dd V(s)}{\dd s}=0,\\
&[X(t),V(t)]=[x,v].\nonumber
\end{aligned}
\right.
\end{equation}
The solution is then given by
\begin{equation}
[X(s),V(s)]=[X(s;t,x,v),V(s;t,x,v)]=[x-(t-s)v^1,v].\nonumber
\end{equation}
We then define the backward exit time $t_{\mathbf{b}}(x,v)$ to be the last moment at which the backward characteristic line $X(-\tau;0,x,v)$ remains in $(0,1)$, that is:
$$t_{\mathbf{b}}(x,v)=\sup \{\tau\geq 0: x-\tau v^1\in (0,1)\}.
$$
We also define $x_{\mathbf{b}}(x,v)=x-t_{\mathbf{b}}v^1\in \{0,1\}.$ Clearly, $v\cdot n(x_{\mathbf{b}}(x,v))\leq 0$.

Let $(x,v)\notin\g_0\cup \g_-$. We set $(t_0,x_0,v_0)=(t,x,v)$. For any $v_{k+1}\in \CV_{k+1}:=\{ v_{k+1}\cdot n(x_k)<0  \}$, the back-time cycle is defined by
\begin{equation}
\left\{\begin{aligned}
X_{cl}(s;t,x,v)&=\sum_{k}\Fi_{[t_{k+1},t_{k})}(s)(x_k-v_{k}^1(t_k-s)),\\[1.5mm]
V_{cl}(s;t,x,v)&=\sum_{k}\Fi_{[t_{k+1},t_{k})}(s)v_k,\label{C.3.3}
\end{aligned}\right.
\end{equation}
with
\begin{equation}
({t}_{k+1},{x}_{k+1},v_{k+1})
=({t}_{k}-{t}_{\mathbf{b}}({x}_{k},v_{k}), {x}_{\mathbf{b}}({x}_{k},v_{k}),v_{k+1}).\nonumber
\end{equation}
We also define the iterated integral
\begin{equation}
\int_{\Pi_{j=1}^{k-1}}\Pi_{j=1}^{k-1}\dd\sigma_j:=\int_{\CV_1}\cdots\left\{\int_{\CV_{k-1}}\dd\sigma_{k-1}\right\}\cdots\dd \sigma_1,\nonumber
\end{equation}
where $\dd\sigma_j:=\sqrt{2\pi}\mu(v_j)|v_{j,1}|\dd v_j$ are the probability measures.

Along the back-cycle \eqref{C.3.3}, we can represent the solution $h$ to the linear equation \eqref{C.2.3.24} in a mild formulation for the $L^\infty$ estimate. Precisely, we have the following lemma. The proof is omitted for brevity as it is similar to that in \cite{EGKM}.
\begin{lemma}[Mild formulation for $h$]
	For any $t>0$ and $(x,v)\in (0,1)\times \mathbb{R}^3\setminus(\gamma_0\cup\gamma_-)$,
	\begin{align}\label{C.3.4}
	h=\sum_{i=1,2,3}J_i+\Fi_{\{t_1>0\}}\sum_{i=4}^{11}J_i,
	\end{align}
	with
	\begin{equation}
	\begin{aligned}
	J_1&=\Fi_{\{t_1\leq  0\}}e^{-\vep^{-1}\nu(v)t}h(x-v_1t,v),\\
	J_2+J_3&=\vep^{-1}\int_{\max\{{t}_1,0\}}^t e^{-\vep^{-1}\nu(v)(t-s)}\left(K_wh+J \right)(x-(t-s)v^1,v)\dd s,\\
	J_4&=e^{-\vep^{-1}\nu(v)(t-t_1)}  q(x_1,v),\\
    J_5&=\frac{e^{-\vep^{-1}\nu(v)(t-t_1)}}{\tilde{w}(v)} \int_{\Pi _{j=1}^{k-1}\mathcal{V}_{j}}
	\sum_{l=1}^{k-2} \Fi_{\{t_{l+1}>0\}}
	q(x_{l+1},v_l)\dd \Sigma_{l}({t}_{l+1}),\\
	J_6&=\frac{e^{-\vep^{-1}\nu(v)(t-t_1)}}{\tilde{w}(v)} \int_{\Pi _{j=1}^{k-1}{\mathcal{V}}_{j}} \sum_{l=1}^{k-1} \Fi_{\{{t}_{l+1}\leq 0<{t}_l\}} h(x_l-v_l^1t_l,v_l) \dd\Sigma_{l}(0),\\
	J_7+J_8&=\frac{e^{-\vep^{-1}\nu(v)(t-t_1)}}{\vep\tilde{w}(v)} \int_{\Pi_{j=1}^{k-1}{\mathcal{V}}_{j}} \sum_{l=1}^{k-1}\int_{0}^{{t}_l} \Fi_{\{{t}_{l+1}\leq 0<t_{l}\}}[K_wh+J](s,x_l-v_l^1(t_l-s),v_l) \dd\Sigma_l(s)\dd s,\\
J_9+J_{10}&=\frac{e^{-\vep^{-1}\nu(v)(t-t_1)}}{\vep\tilde{w}(v)} \int_{\Pi_{j=1}^{k-1}{\mathcal{V}}_{j}} \sum_{l=1}^{k-1}\int_{{t}_{l+1}}^{{t}_l} \Fi_{\{{t}_{l+1}>0\}}[K_wh+J](s,x_l-v_l^1(t_l-s),v_l) \dd\Sigma_l(\s)\dd s,\\
J_{11}&=\frac{e^{-\vep^{-1}\nu(v)(t-t_1)}}{\tilde{w}(v)} \int_{\Pi _{j=1}^{k-1}{\mathcal{V}}_{j}}  \Fi_{\{{t}_{k}>0\}} h({x}_k,v_{k-1}) \dd\Sigma_{k-1}({t}_k),\nonumber
	\end{aligned}
	\end{equation}
where we have denoted
	\begin{align}
	\dd\Sigma_l(s) = \big\{\Pi_{j=l+1}^{k-1}\dd{\sigma}_j\big\}\cdot \big\{\tilde{w}(v_l) e^{-\vep^{-1}\nu(v_l)(t_l-s)} \dd{\sigma}_l\big\}\cdot \big\{\Pi_{j=1}^{l-1} e^{-\vep^{-1}\nu(v_{j})(t_j-t_{j+1})} \dd{\sigma}_j\big\}.\nonumber
	\end{align}
\end{lemma}
\begin{lemma}[cf. \cite{EGKM}]
For $T_0$ sufficiently large, there exists constants $C_1$ and $C_2$ independent of $T_0$, such that for $k=C_1T_0^{\f54}$ and $(x,v)\in (0,1)\times \mathbb{R}_3$,  it holds that
\begin{align}\label{C.3.5}
\int_{\Pi_{j=1}^{k-1}\CV_j}\mathbf{1}_{\{t_k>0\}}\Pi_{j=1}^{k-1}\dd \Sigma_{k-1}(t_k)\leq \left(\frac{1}{2}\right)^{C_2T_0^{\f54}}.	
\end{align}
\end{lemma}	

\begin{lemma}[Weighted $L^\infty$ estimate] \label{lem3.3}
For any $p\in [2,\infty),$ it holds that	
	\begin{align}\label{C.a.1}
	\|h\|_{L^\infty}\leq& C\vep^{-\f1p}\left\|\f{\P_\CM F_R}{\sqrt{\CM}}\right\|_{L^p}+C\vep^{-\f12}\left\|\f{(I-\P_\CM)F_R}{\sqrt{\CM}}\right\|_{L^2}+C|\pa_x\theta_{NS}|_{L^\infty_x}\|h\|_{L^\infty}\nonumber\\
	&+C\left\|\f{\nu^{-1}wg}{\sqrt{\mu}}\right\|_{L^\infty}+C\left|\f{wr}{\sqrt{\mu}}\right|_{L^\infty_-}.
	\end{align}
Here, the constant $C$ is independent of $\vep$.
\end{lemma}	
{\bf Proof:} Take $k=C_1T_0^{\f54}$ so that \eqref{C.3.5} holds. Recall the mild formulation \eqref{C.3.4}. We estimate $J_1-J_{11}$ term by term. Firstly,  we have
\begin{align}\label{C.3.6}
|J_1|\leq Ce^{-\vep^{-1}\nu_0t}\|h\|_{L^\infty}.
\end{align}
For those terms involving $J$ and $r$, notice that
\begin{align}\label{C.3.6-1}
\frac{1}{\tilde{w}(v)}\leq Cw(v)e^{-\frac{|v|^2}{4}}\leq Ce^{-\frac{|v|^2}{8}}.
\end{align}
Thus,  it holds that
\begin{align}
&\int_{\Pi_{j=1}^{k-1}\CV_j}{\bf 1}_{\{t_{l+1}>0\}}\tilde{w}(v_l)\Pi_{j=1}^{k-1}\dd \sigma_j\leq C, ~\text{for }1\leq l\leq k-1, \label{C.3.7}\\
&\int_{\Pi_{j=1}^{k-1}\CV_j}\sum_{l=1}^{k-1}{\bf 1}_{\{t_{l+1}\leq 0<t_l\}}\tilde{w}(v_l)\Pi_{j=1}^{k-1}\dd \sigma_j\leq Ck.\label{C.3.8}
\end{align}
Then by using $$\vep^{-1}\int_{s_1}^{s_2}e^{-\vep^{-1}\nu(v)(s_2-s)}\nu(v)\dd s\leq 1,~\text{for any } s_1<s_2,$$ we can deduce from \eqref{C.3.7} and \eqref{C.3.8} that
\begin{align}
|J_3|+|J_8|+|J_{10}|&\leq Ck\|\nu^{-1}J\|_{L^\infty},\label{C.3.9}\\
|J_4|+|J_5|&\leq Ck|q|_{L^\infty_-}.\label{C.3.10}
\end{align}
By \eqref{C.3.6-1} and \eqref{C.3.8}, we obtain
\begin{align}\label{C.3.12}
|J_6|\leq Cke^{-\vep^{-1}\nu_0t}\|h\|_{L^\infty}.
\end{align}
For $J_{11}$, it follows from \eqref{C.3.5} and \eqref{C.3.6-1} that
\begin{align}\label{C.3.11}
|J_{11}|\leq C\left(\f12\right)^{C_2T_0^{\f54}}\|h\|_{L^\infty}.
\end{align}
For $J_7$, it holds that
\begin{align}
|J_7|\leq& C\vep^{-1}\sum_{l=1}^{k-1}\int_{\Pi_{j=1}^{l-1}\CV_j}\dd \sigma_{j-1}\cdots\dd\sigma_1\int_0^{t_l}e^{-\vep^{-1}\nu_0(t-s)}\dd s\nonumber\\
&\times\int_{\CV_l}\int_{\mathbb{R}_3}{\bf 1}_{\{t_{l+1}\leq 0<t_l\}}\tilde{w}(v_l)|k_w(v_l,v')|h(x_l-v_l^1(t_l-s),v')|\dd v'\dd \sigma_l.\nonumber\\
=&C\vep^{-1}\sum_{l=1}^{k-1}\int_{\Pi_{j=1}^{l-1}\CV_j}\dd \sigma_{j-1}\cdots\dd\sigma_1\int_{t_l-\f{\vep}{N}}^{t_l}\dd s\int_{\CV_l}\int_{\mathbb{R}^3}(\cdots)\dd v'\dd \sigma_l\nonumber\\
&+C\vep^{-1}\sum_{l=1}^{k-1}\int_{\Pi_{j=1}^{l-1}\CV_j}\dd \sigma_{j-1}\cdots\dd\sigma_1\int_{0}^{t_l-\frac{\vep}{N}}\dd s\int_{\CV_l\cap\{|v_l|\geq N \}}\int_{\mathbb{R}^3}(\cdots)\dd v'\dd \sigma_l\nonumber\\
&+C\vep^{-1}\sum_{l=1}^{k-1}\int_{\Pi_{j=1}^{l-1}\CV_j}\dd \sigma_{j-1}\cdots\dd\sigma_1\int_{0}^{t_l-\frac{\vep}{N}}\dd s\int_{\CV_l\cap\{|v_l|\leq N \}}\int_{\{|v'|\geq 2N\}}(\cdots)\dd v'\dd \sigma_l\nonumber\\
&+C\vep^{-1}\sum_{l=1}^{k-1}\int_{\Pi_{j=1}^{l-1}\CV_j}\dd \sigma_{j-1}\cdots\dd\sigma_1\int_{0}^{t_l-\frac{\vep}{N}}\dd s\int_{\CV_l\cap\{|v_l|\leq N \}}\int_{\{|v'|\leq 2N\}}(\cdots)\dd v'\dd \sigma_l\nonumber\\
:=&\sum_{l=1}^{k-1}J_{71l}+J_{72l}+J_{73l}+J_{74l}.\label{C.3.13}
\end{align}
For each term $J_{71l}$, we have
\begin{align}
|J_{71l}|\leq \frac{C}{N}\|h\|_{L^\infty}.\label{C.3.14}
\end{align}
For $J_{72l}$, one  can deduce the following bound
\begin{align}
|J_{72l}|\leq& \f{C}{N}\vep^{-1}\|h\|_{L^\infty}\int_{\Pi_{j=1}^{l-1}\CV_j}\dd\sigma_{j-1}\cdots\dd\sigma_1\nonumber\\
&\times\int_0^{t_l-\f{\vep}{N}}e^{-\vep^{-1}\nu_0(t-s)}\dd s\int_{\CV_l\cap\{|v_l|\geq N \}}e^{-\frac{1}{8}|v_l|^2}\dd v_l\nonumber\\
\leq& Ce^{-\f{N^2}{16}}\|h\|_{L^\infty}.\label{C.3.15}
\end{align}
For $J_{73l}$, we have  $|v'-v_l|\geq N$. Then by \eqref{k2}, we can obtain
\begin{align}\label{C.3.16}
|J_{73l}|\leq& C\vep^{-1}\|h\|_{L^\infty}\int_{\Pi_{j=1}^{l-1}\CV_j}\dd\sigma_{j-1}\cdots\dd\sigma_1\nonumber\\
&\times\int_0^{t_l-\f{\vep}{N}}e^{-\vep^{-1}\nu_0(t-s)}\dd s\int_{\CV_l\cap\{|v_l|\leq N \}}e^{-\frac{1}{8}|v_l|^2}\dd v_l\int_{\mathbb{R}^3}|k_w(v_l,v')|e^{\frac{|v-v'|^2}{64}}e^{-\f{N^2}{64}}\dd v'\nonumber\\
\leq& Ce^{-\f{N^2}{64}}\|h\|_{L^\infty}.
\end{align}
To estimate $J_{74l}$, we introduce the smooth approximate kernel $k_{N}(v_l,v')$ which has a compact support such that
$$\sup_{|v_l|\leq N}\int_{\{|v'|\leq 2N\}}|k_{w}(v_l,v')-k_{N}(v_l,v')|\dd v'\leq \frac{C}{N}.
$$
Then it holds that
\begin{align}\label{C.3.17}
|J_{74l}|\leq& C\vep^{-1}\sum_{l=1}^{k-1}\int_{\Pi_{j=1}^{l-1}\CV_j}\dd \sigma_{j-1}\cdots\dd\sigma_1\int_0^{t_l-\f{\vep}{N}}e^{-\vep^{-1}\nu_0(t-s)}\dd s\nonumber\\
&\times\int_{\CV_l\cap\{|v_l|\leq N\}}\int_{\{|v'|\leq 2N\}}{\bf 1}_{\{t_{l+1}\leq 0<t_l\}}|k_N(v_l,v')h(x_l-v_l^1(t_l-s),v')|\dd v'\dd v_l+\frac{C}{N}\|h\|_{L^\infty}\nonumber\\
\leq &C\vep^{-1}\sum_{l=1}^{k-1}\int_{\Pi_{j=1}^{l-1}\CV_j}\dd \sigma_{j-1}\cdots\dd\sigma_1\int_0^{t_l-\frac{\vep}{N}}e^{-\vep^{-1}\nu_0(t-s)}\dd s\nonumber\\
&\times\int_{\CV_l\cap\{|v_l|\leq N\}}\int_{\{|v'|\leq 2N\}}{\bf 1}_{\{t_{l+1}\leq 0<t_l\}}\left|\f{F_R}{\sqrt{\CM}}(x_l-v_l^1(t_l-s),v')\right|\dd v'\dd v_l+\frac{C}{N}\|h\|_{L^\infty},
\end{align}
where we have used $|k_N(v_l,v')|\leq C_N$. To estimate the first term on the right hand side of \eqref{C.3.17}, we  decompose $F_R=\P_{\CM}F_R+(I-\P_{\CM})F_R$. It then follows from H\"{o}lder inequality that
\begin{align}\label{C.3.18}
&\int_{\CV_l\cap\{|v_l|\leq N\}}\int_{\{|v'|\leq 2N\}}{\bf 1}_{\{t_{l+1}\leq 0<t_l\}}\left|\f{F_R}{\sqrt{\CM}}(x_l-v_l^1(t_l-s),v')\right|\dd v'\dd v_l\nonumber\\
&\quad\leq\left( \int_{\CV_l\cap\{|v_l|\leq N\}}\int_{\{|v'|\leq 2N\}}{\bf 1}_{\{t_{l+1}\leq 0<t_l\}}\left|\f{\P_\CM F_R}{\sqrt{\CM}}(x_l-v_l^1(t_l-s),v')\right|^p\dd v'\dd v_l\right)^{\f1p}\nonumber\\
&\qquad+\left( \int_{\CV_l\cap\{|v_l|\leq N\}}\int_{\{|v'|\leq 2N\}}{\bf 1}_{\{t_{l+1}\leq 0<t_l\}}\left|\f{(I-\P_\CM) F_R}{\sqrt{\CM}}(x_l-v_l^1(t_l-s),v')\right|^2\dd v'\dd v_l\right)^{\f12}.
\end{align}
Set $y_l:=x_l-v_{l}^1(t_l-s)\in (0,1)$ for $s\in [0,t_l-\frac{\vep}{N}].$
By making change of variable $v_l^1\rightarrow y_l$ in the two terms on the right hand side of \eqref{C.3.18} and by noting that the Jacobian
$$
\left|\f{\pa y_l}{\pa v_l^1}\right|=(t_l-s)\geq \frac{\vep}{N},
$$
we deduce that
\begin{align}
&\int_{\CV_l\cap\{|v_l|\leq N\}}\int_{\{|v'|\leq 2N\}}{\bf 1}_{\{t_{l+1}\leq 0<t_l\}}\left|\f{F_R}{\sqrt{\CM}}(x_l-v_l^1(t_l-s),v')\right|\dd v'\dd v_l\nonumber\\
&\qquad\leq C_N\vep^{-\f1p}\left\|\f{\P_\CM F_R}{\sqrt{\CM}}\right\|_{L^p}+C_N\vep^{-\f12}\left\|\f{(I-\P_\CM)F_R}{\sqrt{\CM}}\right\|_{L^2}.\label{C.3.18-1}
\end{align}
By substituting this into \eqref{C.3.17} and combining with \eqref{C.3.13}-\eqref{C.3.16}, we obtain
\begin{align}\label{C.3.19}
|J_7|\leq \frac{Ck}{N}\|h\|_{L^\infty}+C_{N}k\vep^{-\f1p}\left\|\f{\P_\CM F_R}{\sqrt{\CM}}\right\|_{L^p}+C_{N}k\vep^{-\f12}\left\|\f{(I-\P_\CM)F_R}{\sqrt{\CM}}\right\|_{L^2}.
\end{align}
By using the same argument, we have
\begin{align}\label{C.3.20}
|J_9|\leq \frac{Ck}{N}\|h\|_{L^\infty}+C_{N}k\vep^{-\f1p}\left\|\f{\P_\CM F_R}{\sqrt{\CM}}\right\|_{L^p}+C_{N}k\vep^{-\f12}\left\|\f{(I-\P_\CM)F_R}{\sqrt{\CM}}\right\|_{L^2}.
\end{align}
Putting \eqref{C.3.6}, \eqref{C.3.9}-\eqref{C.3.11}, \eqref{C.3.19} and \eqref{C.3.20} together yields the following pointwise estimate
\begin{align}\label{C.3.21}
|h(x,v)|\leq\vep^{-1} \int_{\max\{t_1,0\}}^te^{-\vep^{-1}\nu_0(t-s)}\int_{\mathbb{R}^3}|k_w(v,u)h(x-(t-s)v^1,u)|\dd u\dd s+A(t),
\end{align}
where 
\begin{align}\label{C.3.22}
A(t)=&CT_0^{\f54}\left\{e^{-\vep^{-1}\nu_0t}+\left(\f12\right)^{C_2T_0^{\f54}}+\f1N \right\}\|h\|_{L^\infty}+CT_0^{\f54}\{\|\nu^{-1}J\|_{L^\infty}+|q|_{L^\infty_-}\}\nonumber\\
&+C_{N}T_0^{\f54}\left\{\vep^{-\f1p}\left\|\f{\P_\CM F_R}{\sqrt{\CM}}\right\|_{L^p}+\vep^{-\f12}\left\|\f{(I-\P_\CM)F_R}{\sqrt{\CM}}\right\|_{L^2}\right\}.
\end{align}
Now we denote $y=x-(t-s)v^1\in (0,1)$ and $t_1'=t_1(s,y,u)$ for any $s\in (\max\{t_1,0\},t)$. Then iterating \eqref{C.3.21} once gives
\begin{align}\label{C.3.23}
|h(x,v)|\leq & A(t)+\vep^{-1}\int_{\max\{t_1,0\}}^te^{-\vep^{-1}\nu_0(t-s)}A(s)\dd s\int_{\mathbb{R}^3}|k_w(v,u)|\dd u+B\nonumber\\
\leq& B+CT_0^{\f54}\left\{e^{-\vep^{-1}\nu_0t}(1+t)+\left(\f12\right)^{C_2T_0^{\f54}}+\f1N \right\}\|h\|_{L^\infty}+CT_0^{\f54}\{\|\nu^{-1}J\|_{L^\infty}+|q|_{L^\infty_-}\}\nonumber\\
&+C_{N}T_0^{\f54}\left\{\vep^{-\f1p}\left\|\f{\P_\CM F_R}{\sqrt{\CM}}\right\|_{L^p}+\vep^{-\f12}\left\|\f{(I-\P_\CM)F_R}{\sqrt{\CM}}\right\|_{L^2}\right\},
\end{align}
where
\begin{align}
B:=&\vep^{-2}\int_0^t\dd s\int_0^s e^{-\vep^{-1}\nu_0(t-\tau)}\dd \tau \int_{\mathbb{R}^3}\int_{\mathbb{R}^3}|k_w(v,u)k_w(u,u')|\nonumber\\
&\times{\bf 1}_{\{\max\{t_1,0\}<s<t \}} {\bf 1}_{\{\max\{t_1',0\}<\tau<s \}}|h(y-(s-\tau)u^1,u')|\dd u\dd u'.\nonumber
\end{align}
We estimate $B$  by considering the following cases.\\
{\it Case 1.} For $|v|\geq N$, by \eqref{k1}, we have
\begin{align}\label{C.2.24}
B\leq C(1+|v|)^{-2}\|h\|_{L^\infty}\leq \f{C}{N^2}\|h\|_{L^\infty}.
\end{align}
{\it Case 2.} For $|v|\leq N, |u|\geq 2N$ or $|v|\leq N, |u|\leq 2N, |u'|\geq 3N$, by using \eqref{k2}, we
deduce
\begin{align}\label{C.2.25}
&\vep^{-2}\int_0^t\dd s\int_0^s e^{-\vep^{-1}\nu_0(t-\tau)}\dd \tau\left\{ \int_{\{|u|\geq 2N\}}+\int_{\{|u|\leq 2N,|u'|\geq 3N\}}\right\}(\cdots)\dd u\dd u'\nonumber\\
&\quad\leq e^{-\frac{N^2}{32}}\|h\|_{L^\infty}\left\{ \int_{\{|u|\geq 2N\}}|k_w(v,u)|e^{\f{|v-u|^2}{32}}\dd u+\int_{\{|u|\leq  2N, |u'|\geq 3N\}}|k_w(v,u)k_w(u,u')|e^{\f{|u'-u|^2}{32}}\dd u\right\}\nonumber\\
&\quad\leq e^{-\frac{N^2}{32}}\|h\|_{L^\infty}.
\end{align}	
{\it Case 3.}  For  $|v|\leq N, |u|\leq 2N, |u'|\leq 3N$, we have
\begin{align}\label{C.2.26}
&\vep^{-2}\int_0^t\dd s\int_0^s e^{-\vep^{-1}\nu_0(t-\tau)}\dd \tau\int_{\{|u|\leq 2N,|u'|\leq 3N\}}(\cdots)\dd u\dd u'\nonumber\\
&\qquad\leq  \frac{C}{N}\|h\|_{L^\infty}+\vep^{-2}\int_0^t\dd s\int_0^{s-\frac{\vep}{N}} e^{-\vep^{-1}\nu_0(t-\tau)}\dd \tau\int_{\{|u|\leq 2N,|u'|\leq 3N\}}|k_N(v,u)k_N(u,u')|\nonumber\\
&\qquad\quad\times{\bf 1}_{\{\max\{t_1,0\}<s<t \}} {\bf 1}_{\{\max\{t_1',0\}<\tau<s \}}|h(y-(s-\tau)u^1,u')|\dd u\dd u'.\nonumber\\
&\qquad\leq \frac{C}{N}\|h\|_{L^\infty}+C_N\vep^{-2}\int_0^t{\bf 1}_{\{\max\{t_1,0\}<s<t \}}\dd s\int_0^{s-\frac{\vep}{N}} e^{-\vep^{-1}\nu_0(t-\tau)}\dd \tau\nonumber\\
&\qquad\quad\times \int_{\{|u|\leq 2N,|u'|\leq 3N\}}{\bf 1}_{\{\max\{t_1',0\}<\tau<s \}}\left|\f{F_R}{\sqrt{\CM}}(y-(s-\tau)u^1,u')\right|\dd u\dd u',
\end{align}
where we have used the boundedness of smooth approximate kernel $k_N(v,u).$ For the last term on the right hand side of \eqref{C.2.26}, we decompose $F_R=\P_\CM F_R+(I-\P_{\CM})F_R$. Then similar to \eqref{C.3.18-1}, we obtain
\begin{align}\label{C.2.26-1}
&\int_{\{|u|\leq 2N,|u'|\leq 3N\}}{\bf 1}_{\{\max\{t_1',0\}<\tau<s \}}\left|\f{F_R}{\sqrt{\CM}}(y-(s-\tau)u^1,u')\right|\dd u\dd u'\nonumber\\
&\qquad \leq C \vep^{-\f1p}\left\|\f{\P_\CM F_R}{\sqrt{\CM}}\right\|_{L^p}+\vep^{-\f12}\left\|\f{(I-\P_\CM)F_R}{\sqrt{\CM}}\right\|_{L^2}.
\end{align}
By combining \eqref{C.2.24}-\eqref{C.2.26-1}, we have
\begin{align}
B\leq \frac{C}{N}\|h\|_{L^\infty}+\leq C \vep^{-\f1p}\left\|\f{\P_\CM F_R}{\sqrt{\CM}}\right\|_{L^p}+\vep^{-\f12}\left\|\f{(I-\P_\CM)F_R}{\sqrt{\CM}}\right\|_{L^2}.\label{C.2.27}
\end{align}
By putting \eqref{C.2.27} into \eqref{C.3.23}, we get
\begin{align}\label{C.3.28}
\|h\|_{L^\infty}
\leq& CT_0^{\f54}\left\{e^{-\vep^{-1}\nu_0t}(1+t)+\left(\f12\right)^{C_2T_0^{\f54}}+\f1N \right\}\|h\|_{L^\infty}+CT_0^{\f54}\{\|\nu^{-1}J\|_{L^\infty}+|q|_{L^\infty_-}\}\nonumber\\
&+C_{N}T_0^{\f54}\left\{\vep^{-\f1p}\left\|\f{\P_\CM F_R}{\sqrt{\CM}}\right\|_{L^p}+\vep^{-\f12}\left\|\f{(I-\P_\CM)F_R}{\sqrt{\CM}}\right\|_{L^2}\right\}.
\end{align}
Now we choose  $N=2CT_0^{\f54}$ and fix $t=T_0$ large enough such that
\begin{align}
CT_0^{\f54}\left\{e^{-\vep^{-1}\nu_0t}(1+t)+\left(\f12\right)^{C_2T_0^{\f54}}+\f1N \right\}\leq \f34.\nonumber
\end{align}
Then it follows from \eqref{C.3.28} that
\begin{align}\label{C.3.29}
\|h\|_{L^\infty}\leq C\vep^{-\f1p}\left\|\f{\P_\CM F_R}{\sqrt{\CM}}\right\|_{L^p}+C\vep^{-\f12}\left\|\f{(I-\P_\CM)F_R}{\sqrt{\CM}}\right\|_{L^2}+C\|\nu^{-1}J\|_{L^\infty}+|q|_{L^\infty_-}.
\end{align}
As for the last two terms in \eqref{C.3.29}, by \eqref{g1}, \eqref{C.3.1} and \eqref{C.3.2}, we have
\begin{align}
\|\nu^{-1}J\|_{L^\infty}&\leq C\left\|\frac{w(\CM-\mu)}{\sqrt{\mu}}\right\|_{L^\infty}\|h\|_{L^\infty}+ C\left\|\f{\nu^{-1}wg}{\sqrt{\mu}}\right\|_{L^\infty}\nonumber\\
&\leq C|\pa_x\theta_{NS}|_{L^\infty_x}\|h\|_{L^\infty}+C\left\|\f{\nu^{-1}wg}{\sqrt{\mu}}\right\|_{L^\infty},\nonumber
\end{align}
and
\begin{align}
|q|_{L^\infty_-}\leq C|\pa_x\theta_{NS}|_{L^\infty}\|h\|_{L^\infty}+C\left|\f{wr}{\sqrt{\mu}}\right|_{L^\infty_-}.\nonumber
\end{align}
Combining these two estimates with \eqref{C.3.29} yields \eqref{C.a.1}. The proof of Lemma \ref{lem3.3} is completed.\qed

We are now  ready to prove the main result in this section.
	
{\bf Proof of Proposition \ref{propC.2.1}:}
The existence part can be established by the same iteration procedure as  in \cite{EGKM,EGKM-hy} and we omit it for brevity. Here, we only show the a priori estimate \eqref{C.2.3.4}. A suitable combination of \eqref{leme1}, \eqref{C.2.3.23} and \eqref{C.a.1} yields
\begin{align}
&\left\|\f{wF_R}{\sqrt{\mu}}\right\|_{L^\infty}+\vep^{-\f1p}\left\|\f{\P_\CM F_R}{\sqrt{\CM}}\right\|_{L^p}+\vep^{-1-\f1p}\left\|\f{\nu^{\f12}(I-\P_\CM)F_R}{\sqrt{\CM}}\right\|_{L^2}+\vep^{-\f12-\f1p}\left|\f{\nu^{\f12}(I-\P_\g)F_R}{\sqrt{\CM}}\right|_{L^2_+}\nonumber\\
&\qquad\leq  C|\pa_x\theta_{NS}|_{L^\infty}\left( \left\|\f{wF_R}{\sqrt{\mu}}\right\|_{L^\infty}+\vep^{-\f1p}\left\|\f{\P_\CM F_R}{\sqrt{\CM}}\right\|_{L^p}+\vep^{-1-\f1p}\left\|\f{\nu^{\f12}(I-\P_\CM)F_R}{\sqrt{\CM}}\right\|_{L^2} \right)\nonumber\\
&\qquad\quad+C\vep^{-\f12}\left\|\f{\nu^{\f12}(I-\P_\CM)F_R}{\sqrt{\CM}}\right\|_{L^2}+C\left(\vep^{\f{p+2}{p(p-2)}}+\vep^{-1-p}e^{-\f{\tau_0}{\vep}}\right)\left\|\f{wF_R}{\sqrt{\mu}}\right\|_{L^\infty}
\nonumber\\
&\quad\qquad+
C\vep^{-\f1p}\left|\f{\nu^{\f12}(I-\P_\g)F_R}{\sqrt{\CM}}\right|_{L^2_+}+C\vep^{-2-\f1p}\left\|\f{\P_{\CM}g}{\sqrt{\CM}}\right\|_{L^2}+C\vep^{-1-\f1p}\left\|\f{\nu^{-\f12}(I-\P_\CM)g}{\sqrt{\CM}}\right\|_{L^2}\nonumber\\
&\qquad\quad+C\vep^{-\f12-\f1p}\left|\f{r}{\sqrt{\CM}}\right|_{L^2_{-}}+C\left\|\f{\nu^{-1}wg}{\sqrt{\mu}}\right\|_{L^\infty}+C\left|\f{wr}{\sqrt{\mu}}\right|_{L^\infty_-}.\label{C.3.30}
\end{align}
Taking both $|\theta_0-\theta_1|$ and $
\vep$ suitably small such that
\begin{align}
C|\pa_x\theta_{NS}|_{L^\infty}\leq \frac{1}{4},~ ~C\left(\vep^{\f{p+2}{p(p-2)}}+\vep^{-1-p}e^{-\f{\tau_0}{\vep}}\right)\leq \f14, \nonumber
\end{align}
we can absorb all of the $F_R$ terms on the right hand side of \eqref{C.3.30} by the left hand side. Thus, the a priori estimate \eqref{C.2.3.4} follows and the proof of Proposition \ref{propC.2.1} is completed. \qed
	\section{Justification of the expansion \label{sec5}}
	In this section, we will solve the remainder system \eqref{C.2.2.3} with boundary condition \eqref{C.2.2.4} and then give the proof of Theorem \ref{thmC.1.1}.
	\begin{lemma}\label{lmC.2.4.1}
	Let $w$ be the weight function defined in \eqref{w}. Then it holds 
		\begin{align}
		&\left\|\f{\nu^{-\f12}L_{as}F_R}{\sqrt{\CM}}\right\|_{L^2} \leq C|\pa_x\theta_{NS}|_{L^\infty_x}\cdot\left\|\f{\nu^{\f12}F_R}{\sqrt{\CM}}\right\|_{L^2},\label{C.2.4.1-4}\\ &\left\|\f{w\nu^{-1}L_{as}F_R}{\sqrt{\mu}}\right\|_{L^\infty}\leq  C|\pa_x\theta_{NS}|_{L^\infty_x}\cdot\left\|\f{wF_R}{\sqrt{\mu}}\right\|_{L^\infty},\label{C.2.4.1-5}\\
		&\left\|\f{\nu^{-\f12}Q(F_R,F_R)}{\sqrt{\CM}}\right\|_{L^2} +\left\|\f{w\nu^{-1}Q(F_R,F_R)}{\sqrt{\mu}}\right\|_{L^\infty} \leq C\left\|\f{wF_R}{\sqrt{\mu}}\right\|^2_{L^\infty}.\label{C.2.4.1-6}
		\end{align}
	\end{lemma}
{\bf Proof:} By \eqref{g1}, \eqref{6.2} and \eqref{6.3}, we have
\begin{align}
\left\|\f{\nu^{-\f12}L_{as}F_R}{\sqrt{\CM}}\right\|_{L^2} &\leq C\left\|\f{\nu^{\f12}F_R}{\sqrt{\CM}}\right\|_{L^2}\cdot \left( \left\|\f{\nu^{\f12}(G,\vep F_2)}{\CM}\right\|_{L^2} +\left\|\f{\nu^{\f12}(\mathfrak{B}_0,\mathfrak{B}_1)}{\sqrt{\CM}}\right\|_{L^2} \right)\nonumber\\
&\leq C|\pa_x\theta_{NS}|_{L^\infty_x}\left\|\f{\nu^{\f12}F_R}{\sqrt{\CM}}\right\|_{L^2},\nonumber
\end{align}
which is \eqref{C.2.4.1-4}.
\eqref{C.2.4.1-5} can be obtained similarly. Finally, by \eqref{g2} and the fact that
\begin{align}
\left|\f{\nu^{-\f12}Q(F_R,F_R)}{\sqrt{\CM}}\right|_{L^2_v}\leq \left|w^{-1}\sqrt{\f{\mu}{\CM}}\right|_{L^2_v}\left|\frac{\nu^{-1}wQ(F_R,F_R)}{\sqrt{\mu}}\right|_{L^\infty_v}\leq C \left|\frac{\nu^{-1}wQ(F_R,F_R)}{\sqrt{\mu}}\right|_{L^\infty_v},\nonumber
\end{align}
we can deduce
\eqref{C.2.4.1-6}. Therefore, the proof of Lemma \ref{lmC.2.4.1} is completed.

{\bf Proof of Theorem \ref{thmC.1.1}.} We construct the solution to the remainder system \eqref{C.2.2.3} with boundary condition \eqref{C.2.2.4} via the following iteration scheme:
	\begin{equation}\label{C.2.4.2}
	\left\{
	\begin{aligned}
	&\vep v_1\pa_xF_R^{n+1}+\CL_{\CM}F_R^{n+1}=-\vep \CL_{as}F_R^{n}+\vep^{1+\alpha}Q(F_R^{n},F_R^n)+\vep^{1-\alpha}A_s,\\
	&F_R^{n+1}|_{\g_-}=\P_\gamma F_R^{n+1}+\mathbb{I}_\g F_R^n+\vep^{1-\alpha}r,\\
	&F_R^0\equiv 0,
	\end{aligned}\right.
	\end{equation}
	where
	$$\mathbb{I}_{\g}F_R^n=\left(\mu_w-\frac{1}{2\pi\theta_{NS}^2}e^{-\f{|v|^2}{2\theta_{NS}}}\right)\int_{\{n(x)\cdot v>0\}}F_R^{n}\{n(x)\cdot v\}\dd v.
	$$
	Direct computation yields that
	$$\int_{0}^1\int_{\mathbb{R}^3}\CL_{as}F_R^{n}\dd x\dd v=\int_{0}^1\int_{\mathbb{R}^3}Q(F_R^{n},F_R^n)\dd x\dd v=\int_{0}^1\int_{\mathbb{R}^3}A_s\dd x\dd v=0,
	$$
	$$
	\begin{aligned}
	\int_{\{v_1>0\}}\mathbb{I}_{\g} F_R^n(0,v)v_1\dd v&=\int_{\{v_1>0\}}r(0,v)v_1\dd v=0,\\
	\int_{\{v_1<0\}}\mathbb{I}_{\g} F_R^n(1,v)v_1\dd v&=\int_{\{v_1<0\}}r(1,v)v_1\dd v=0,
	\end{aligned}
	$$
	and
	$$\begin{aligned}
	\left|\frac{\mathbb{I}_\g F_R^n}{\sqrt{\CM}}\right|_{L^2_-}+\left|\frac{w\mathbb{I}_\g F_R^n}{\sqrt{\mu}}\right|_{L^\infty_-}&\leq  C\max\{|\theta_{NS}(0)-\theta_0|,|\theta_{NS}(1)-\theta_1|\}\cdot\left|\f{wF^n_R}{\sqrt{\mu}}\right|_{L^\infty_+}\\
	&\leq C\vep|\pa_{x}\theta_{NS}|_{L^\infty_x}\cdot\left|\f{wF^n_R}{\sqrt{\mu}}\right|_{L^\infty_+}.
	\end{aligned}
	$$
	Therefore, the existence of sequence of solutions $F_R^n$ to the system \eqref{C.2.4.2} follows from Proposition \ref{propC.2.1}. Then applying the estimate \eqref{C.2.3.4} to $F_{R}^{n+1}$, we have
	\begin{align}
	&\left\|\f{wF_R^{n+1}}{\sqrt{\mu}}\right\|_{L^\infty}+\vep^{-\f1p}\left\|\f{\P_{\CM}F_R^{n+1}}{\sqrt{\CM}}\right\|_{L^p}+\vep^{-1-\f1p}\left\|\f{\nu^{\f12}(I-\P_\CM)F_R^{n+1}}{\sqrt{\CM}}\right\|_{L^2}\nonumber\\
	&\quad\lesssim_p \vep^{-1-\f1p}\left[\vep\left\|\f{\nu^{-\f12}\CL_{as}F_R^n}{\sqrt{\CM}}\right\|_{L^2}+\vep^{1+\alpha}\left\|\f{\nu^{-\f12}Q(F_R^n,F_R^n)}{\sqrt{\CM}}\right\|_{L^2}
	+\vep^{1-\alpha}\left\|\f{\nu^{-\f12}A_s}{\sqrt{\CM}}\right\|_{L^2}\right]\nonumber\\
	&\qquad+\vep^{\f12-\f1p-\alpha}\left[\left|\f{r}{\sqrt{\CM}}\right|_{L^2_{-}}+\left|\f{wr}{\sqrt{\mu}}\right|_{L^\infty_-}\right]+\vep^{\f12-\f1p}|\pa_x\theta_{NS}|_{L^\infty_x}\cdot\left|\f{wF_R^n}{\sqrt{\mu}}\right|_{L^\infty_{+}}
	\nonumber\\
	&\qquad+\left[\vep\left\|\f{w\nu^{-1}\CL_{as}F_R^n}{\sqrt{\mu}}\right\|_{L^\infty}+\vep^{1+\alpha}\left\|\f{w\nu^{-1}Q(F_R^n,F_R^n)}{\sqrt{\mu}}\right\|_{L^\infty}
	+\vep^{1-\alpha}\left\|\f{w\nu^{-1}A_s}{\sqrt{\mu}}\right\|_{L^\infty}\right].\nonumber
	\end{align}
	Here, we have used the fact that $\P_\CM A_s=0$. By using bounds \eqref{6.4}, \eqref{6.5} in Lemma \ref{lem2.1} and \eqref{C.2.4.1-4}-\eqref{C.2.4.1-6} in Lemma \ref{lmC.2.4.1}, we can further obtain
	\begin{align}
	&\left\|\f{wF_R^{n+1}}{\sqrt{\mu}}
	\right\|_{L^\infty}+\vep^{-\f1p}\left\|\f{\P_{\CM}F_R^{n+1}}{\sqrt{\CM}}\right\|_{L^p}+\vep^{-1-\f1p}\left\|\f{\nu^{\f12}(I-\P_\CM)F_R^{n+1}}{\sqrt{\CM}}\right\|_{L^2}\nonumber\\
	&\quad\lesssim_p \vep^{-\f1p}|\pa_{x}\theta_{NS}|_{L^\infty}\left[\left\|\frac{\P_\CM F_R^n}{\sqrt{\CM}}\right\|_{L^p}+\left\|\frac{\nu^{\f12}(I-\P_\CM)F_R^n}{\sqrt{\CM}}\right\|_{L^2}\right]
	\nonumber\\
	&\qquad+\vep^{\alpha-\f1p}\left\|\f{wF_R^n}{\sqrt{\mu}}
	\right\|_{L^\infty}^2+\vep^{\f12-\f1p}\left\|\f{wF_R^n}{\sqrt{\mu}}
	\right\|_{L^\infty}+C_p\vep^{\f12-\f1p-\alpha}|\pa_x\theta_{NS}|_{L^\infty_x}.\nonumber
	\end{align}
	We now  fix $0<\alpha<1/2$. Taking $p>2$ suitably large and $|\theta_1-\theta_0|$ suitably small, we have the following estimate
	$$\left\|\f{wF_R^{n+1}}{\sqrt{\mu}}
	\right\|_{L^\infty}+\vep^{-\f1p}\left\|\f{\P_{\CM}F_R^{n+1}}{\sqrt{\CM}}\right\|_{L^p}+\vep^{-1-\f1p}\left\|\f{\nu^{\f12}(I-\P_\CM)F_R^{n+1}}{\sqrt{\CM}}\right\|_{L^2}\leq C|\pa_x\theta_{NS}|_{L^\infty_x},
	$$
	where the constant $C>0$ is independent of $n.$ Moreover, it is straightforward to show $F_R^{n+1}$ is a Cauchy sequence. Hence the solution to the remainder system \eqref{C.2.2.3} with boundary condition \eqref{C.2.2.4} is constructed by taking $n\rightarrow\infty$ and the estimate \eqref{C.1.3.2} follows immediately. Therefore, the proof of Theorem \ref{thmC.1.1} is completed. \qed

\section{Appendix}
The following Lemma summarizes the well-posedness of Milne problem in $L^\infty$ space that was proved in  \cite{BCN,Wu}.
\begin{lemma}\label{lmC.a.1}
	Let $0\leq \varpi<1/4$ and $\beta>3$. Suppose that $$\sup_{v_1>0}\left|(1+|v|^2)^{\f\beta2}e^{\varpi|v|^2}\frac{\mathcal{G}}{\sqrt{\mu}}\right|<\infty.$$ Then there exist a positive constant $\sigma_0>0$ and a smooth function $\mathcal F_{\infty}\in \text{Ker}\CL_{\mu}$, such that \eqref{C.a.2} admits a unique solution $\mathcal{F}$ satisfying
	\begin{align}\label{C.a.3}
		\left|(1+|\cdot|^2)^{\f\beta2}e^{\varpi|\cdot|^2}\frac{\mathcal{F}(y)-\mathcal{F}_{\infty}}{\sqrt{\mu}}\right|_{L^\infty_v}\leq Ce^{-\sigma_0y}\sup_{v_1>0}\left|(1+|v|^2)^{\f\beta2}e^{\varpi|v|^2}\frac{\mathcal{G}}{\sqrt{\mu}}\right|.
	\end{align}
\end{lemma}

\bigskip

{\bf Acknowledgement}\ \ Renjun Duan's research was partially supported by the General Research Fund (Project No.~14301719) from RGC of Hong Kong and a Direct Grant (4053397) from CUHK. Shuangqian Liu's research was supported by grants from the National Natural Science Foundation of China (contracts: 11971201 and 11731008).
Tong Yang's research was supported by a fellowship award from the Research Grants Council of the
 Hong Kong Special Administrative Region, China (Project no. SRF2021-1S01). The  authors would like to thank Professor Kazuo Aoki for introducing the problem as well as pointing out Coron's paper \cite{Coron} in 2018.

\end{document}